\documentclass[12pt,twoside,leqno]{amsart}
\usepackage{amssymb,amsbsy,amsmath,amsfonts,amssymb,amscd,epsfig,
times,graphics,color,xypic}

\definecolor{blue}{cmyk}{1.,1.,0.,0.63}
\definecolor{red}{cmyk}{0.,1.,1.,0.63}
\definecolor{green}{cmyk}{1.,0.,1.,0.63}

\usepackage[latin1]{inputenc}
\newcommand{\B}{\mathbb{B}}
\newcommand{\C}{\mathbb{C}}

\newcommand{\N}{\mathbb{N}}
\newcommand{\R}{\mathbb{R}}

\newcommand{\blue}{\textcolor{blue}}
\newcommand{\green}{\textcolor{green}}
\newcommand{\red}{\textcolor{red}}

\newtheorem{definition}{Definition}
\newtheorem{theorem}{Theorem}
\newtheorem{proposition}{Proposition}
\newtheorem{lemma}{Lemma}
\newtheorem{corollary}{Corollary}

\begin{document}


\title[A Morse-theoretical proof
of the Hartogs extension theorem]{
A Morse-theoretical proof
\\
of the Hartogs extension theorem
}

\author{Jo\"el Merker and Egmont Porten}

\address{
D\'epartement de Math\'ematiques et Applications, UMR 8553
du CNRS, \'Ecole Normale
Sup\'erieure, 45 rue d'Ulm, F-75230 Paris Cedex 05, 
France. \ \
{\it Internet}:
{\tt http://www.cmi.univ-mrs.fr/$\sim$merker/index.html}}

\email{merker@dma.ens.fr}

\address{Department of Engineering, Physics and
Mathematics, Mid Sweden University, Campus Sundsvall,
S-85170 Sundsvall, Sweden
}

\email{Egmont.Porten@miun.se}

\date{\number\year-\number\month-\number\day}

\begin{abstract}
100 years ago exactly, in 1906, Hartogs published a celebrated
extension phenomenon (birth of {\em Several Complex Variables}), 
whose global counterpart was understood later: 
{\it holomorphic functions in a connected neighborhood
$\mathcal{ V} ( \partial \Omega)$ of a connected boundary $\partial
\Omega \Subset \C^n$ {\rm (}$n \geqslant 2${\rm )} do extend
holomorphically and uniquely to the domain $\Omega$}. Martinelli in
the early 1940's and Ehrenpreis in 1961 obtained a rigorous proof,
using a new multidimensional integral kernel or a short $\overline{
\partial}$ argument, but it remained unclear how to derive
a proof using only analytic discs, as did Hurwitz (1897), Hartogs
(1906) and E.E.~Levi (1911) in some special, model cases. In fact,
known attempts ({\it e.g.} Osgood 1929, Brown 1936) struggled for
monodromy against multivaluations, but failed to get the general
global theorem.

Moreover, quite unexpectedly, Forn{\ae}ss in 1998 exhibited a
topologically strange (nonpseudoconvex) domain $\Omega^{\sf F} \subset
\C^2$ that cannot be filled in by holomorphic discs, when one makes
the additional requirement that discs must all lie entirely inside
$\Omega^{ \sf F}$. However, one should point out that the standard,
unrestricted disc method usually allows discs to go outsise the domain
(just think of Levi pseudoconcavity).

Using the method of analytic discs for local extensional steps and
Morse-theoretical tools for the global topological control of
monodromy, we show that the Hartogs extension theorem can be
established in such a way.

\end{abstract}

\maketitle

\vspace{-0.5cm}

\begin{center}
\begin{minipage}[t]{11cm}
\baselineskip =0.35cm
{\scriptsize

\centerline{\bf Table of contents}

\smallskip

{\bf 1.~The Hartogs extension theorem and the method of analytic discs 
\dotfill 1.}

{\bf 2.~Preparation of the boundary and unique extension \dotfill 5.}

{\bf 3.~Quantitative Hartogs-Levi extension by pushing analytic discs 
\dotfill 9.}

{\bf 4.~Filling domains outside balls of decreasing radius \dotfill 13.}

{\bf 5.~Creating domains, merging and suppressing connected components 
\dotfill 20.}

{\bf 6.~The exceptional case $k_\lambda = 1$ \dotfill 30.}

\smallskip

{\footnotesize\tt \hfill \red{[22 colored illustrations]}}

}\end{minipage}
\end{center}

\section*{ \S1.~The Hartogs extension theorem 
\\
and the method of analytic discs}

100 years ago exactly, in 1906, the publication of Hartogs's thesis
(\cite{ ha1906} under the direction of Hurwitz) revealed what is now
considered to be the most striking fact of multidimensional complex
analysis: the automatic, compulsory holomorphic extension of functions
of several complex variables to larger domains, especially for a class
of ``pot-looking'' domains, nowadays called {\it Hartogs figures},
that may be filled in up to their top. Soon after, E.E. Levi \cite{
le1911} applied the Hurwitz-Hartogs argument of Cauchy integration on
complex affine circles moving in the domain (firstly discovered in
\cite{ hu1897}), in order to perform local holomorphic extension
across strictly (Levi) pseudoconcave boundaries. The so-called {\em
method of analytic discs}\, was born, historically.

\medskip\noindent
{\bf Hartogs extension theorem.}
{\it
Let $\Omega \Subset \C^n$ be a bounded domain having {\rm connected}
boundary. \underline{If $n\geqslant 2$}, every function holomorphic in
some {\it connected}\, open neighborhood $\mathcal{ V} ( \partial
\Omega)$ of $\partial \Omega$ extend holomorphically and uniquely
inside $\Omega$, {\it i.e.}{\rm :}}
\[
\forall\,f\in \mathcal{O}
\big(\mathcal{V}(\partial\Omega)\big),
\ \ \ 
\exists\,!\
F\in\mathcal{O}
\big(
\Omega
\cup
\mathcal{V}(\partial\Omega)
\big)
\ \ \
\text{\rm s.t.}
\ \ \
F\big\vert_{
\mathcal{V}(\partial\Omega)}
=f.
\]

Classically, one also presents an alternative formulation, which is
checked to be equivalent\,\,---\,\,think that $K = \Omega \big
\backslash \mathcal{ V} ( \partial \Omega)$.

\medskip\noindent
$\text{\bf Hartogs theorem}^\text{\bf bis}\text{\bf .}$ {\it If
$\Omega \Subset \C^n$ $(n\geqslant 2)$ is a domain and if $K \subset
\Omega$ is any compact such that $\Omega \backslash K$ {\sl
connected}, then $\mathcal{ O} (\Omega \backslash K) = \mathcal{ O}
(\Omega) \big\vert_{ \Omega \backslash K}$}.

\medskip

Already in~\cite{ ha1906} (p.~231), Hartogs stated such a global
theorem in the typical language of those days, without claiming
single-valuedness however\,\,---\,\,something that he consistently
mentions in other places. Later in~\cite{ os1929}, Osgood (who gives
the reference to Hartogs) ``proves'' unique holomorphic extension with
discs, but what is written there is seriously erroneous, even when
applied to a ball. In 1936, well before Milnor (\cite{ mi1963}) had
popularized Morse theory, using topological concepts and a language
which are nowadays difficult to grasp, Brown (\cite{ br1936}) fixed
somehow single-valuedness of the extension\footnote{\, We thank an
anonymous referee for pointing historical incorrections in the
preliminary version of this paper and for providing us with exact
informations.}: discretizing $\Omega \backslash K$ to tame the
topology, he exhausts $\C^n$ by spheres of decreasing radius (as we
will do in this paper), but we believe that his proof still contains
imprecisions, because the subtracting process that we encounter
unavoidably when applying Morse theory with the same
spheres does not appear in~\cite{
br1936}.

Since the 1940's, few complex analysts have seriously thought about
testing the limit of the disc method probably because the motivation
was gone, and in fact, the possible existence of an elementary {\it
rigorous proof}\, of the global Hartogs extension theorem using only a
finite number of Hartogs figures remained a folklore belief; for
instance, in \cite{ st1988}, p.~133, it is just left as an
``exercise''. But to the authors' knowledge, no reliable mathematical
publication shows fully how to perform a rigorous proof of the global
theorem, using only the original Hurwitz-Hartogs-Levi analytic discs
as a tool.

\smallskip
On the other hand, thanks to the contributions 
of Fueter (\cite{ fu1939}), of Martinelli (\cite{ ma1938,
ma1942/43}), of Bochner (\cite{ bo1943}) and of Fichera (\cite{
fi1957}), powerful multidimensional integral kernels were discovered
that provided a complete proof, from the side of Analysis. Soon after,
Ehrenpreis (\cite{ eh1961}) found what is known to be the most concise
proof, based on the vanishing of $\overline{ \partial}$-cohomology
with compact support. This proof was learnt by generations of complex
analysts, thanks to H\"ormander's book~\cite{ ho1966}. Range's {\em
Correction of the Historical Record}~\cite{ ra2002} provides an
excellent account of the very birth of integral formulas in
$\C^n$. Since the 1960's, $\overline{ \partial}$ techniques, $L^2$
methods and integral kernels developed into a vast field of research
in Several Complex Variables, {\it c.f.} \cite{ ho1966, ah1972,
hl1984, he1985, ra1986, cst1994, de1997, japf2000, jp2002, lm2002,
hm2002}.

\smallskip
A decade ago, Forn{\ae}ss~\cite{ fo1998} produced a topologically
strange domain $\Omega^{ \sf F}$ that cannot be filled in by means of
analytic discs, when one makes the additional requirement that discs
must all lie {\it entirely inside}\, the domain. Possibly, one could
interpret this example as a ``defeat'' of geometrical methods.

\smallskip
But in absence of pseudoconvexity, it is much more natural to allow
discs to {\it go outside}\, the domain, because the local E.E.~Levi
extension theorem already needs that. In fact, as remarked by Bedford
in his review~\cite{ be1999} of~\cite{ fo1998}, Hartogs' phenomenon
for Forn{\ae}ss' domain $\Omega^{ \sf F}$ may be shown to hold
straightforwardly by means of the usual, unrestricted disk method.

\smallskip
Furthermore, the study of envelopes of holomorphy ({\it see}\, the
monograph of Jarnicki and Pflug~\cite{ japf2000} for an 
introduction to Riemann domains spread over $\C^n$
and~\cite{ mp2006} for applications in a CR context) shows well how
natural it is to deal with sucessively enlarged (Riemann)
domains. Bishop's constructive approach, especially his famous idea of
gluing discs to real submanifolds, reveals to be adequate in such a
widely open field of research. We hence may hope that, after the very
grounding historical theorem of Hartogs has enjoyed a renewed proof,
geometrical methods will undergo further developments, especially to
devise fine holomorphic extension theorems that are unreachable by
means of contemporary $\overline{ \partial }$ techniques.

\medskip
In this paper, we establish rigorously that the Hartogs extension
theorem can be proved by means of a {\it finite number}\, of
parameterized families of analytic discs (Theorems~2.7 and 5.4).
The discs we use are all (tiny) pieces of complex lines
in $\C^n$. The main difficulty is topological 
and we use the Morse machinery to tame multisheetedness.

\smallskip
At first, we shall replace the boundary $\partial \Omega$ by a
$\mathcal{ C}^\infty$ connected oriented hypersurface $M \Subset \C^n$
($n\geqslant 2$) for which the restriction to $M$ of the Euclidean
norm function $z \mapsto \vert \! \vert z \vert \! \vert$ is a good
Morse function (Lemma~3.3), namely there exist only finitely many
points $\widehat{ p}_\lambda \in M$, $1 \leqslant \lambda \leqslant
\kappa$, with $\vert \! \vert \widehat{ p}_1 \vert \! \vert < \cdots <
\vert \! \vert \widehat{ p}_\kappa \vert \! \vert$ at which $z \mapsto
\vert \! \vert z \vert \! \vert$ restricted to $M$ has vanishing
differential. We also replace $\mathcal{ V} ( \partial \Omega)$ by a
very thin tubular neighborhood $\mathcal{ V}_\delta ( M)$, $0 < \delta
< \! \! < 1$, and $\Omega$ by a domain $\Omega_M \Subset \C^n$ bounded
by $M$. Next, we will introduce a modification of the Hartogs figure,
called a {\sl Levi-Hartogs figure}, which is more appropriate to
produce holomorphic extension from the cut out domains $\big\{ \vert
\! \vert z \vert \! \vert > r \big\} \cap \Omega_M$, where the radius
$r$ will decrease, inductively. Local Levi pseudoconcavity of the
exterior of a ball then enables us to prolong the holomorphic
functions to $\big\{ \vert \! \vert z \vert \! \vert > r - \eta
\big\} \cap \Omega_M$, for some uniform $\eta$ with $0 < \eta < \! \!
< 1$, which depends on the dimension $n \geqslant 2$, on $\delta$, and
on the diameter of $\overline{ \Omega}$. We hence descend stepwise to
lower radii until the domain is fully filled in.

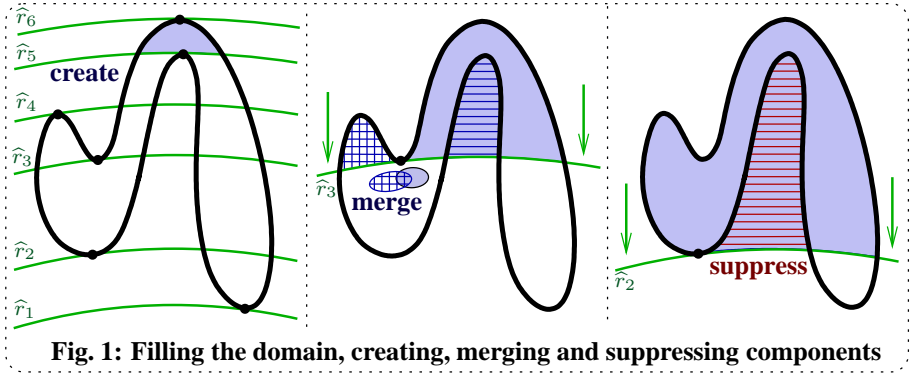
\begin{center}
\input heuristic.pstex_t
\end{center}

However, this naive conclusion fails because
of multivaluations and a crucial three-piece
topological device is required. We begin by filling the top of the
domain, which is simply diffeomorphic to a cut out piece of
ball. Geometrically speaking, Morse points $\widehat{ p}_\lambda$, $1
\leqslant \lambda \leqslant \kappa$, are the only points of $M$ at
which the family of spheres $\big( \big\{ \vert \! \vert z \vert \!
\vert = r \big\}\big)_{ 0 < r < \infty}$ are tangent to $M$. We denote
$\vert \! \vert \widehat{ p}_\lambda \vert \! \vert =: \widehat{
r}_\lambda$ with $\widehat{ r}_1 < \cdots < \widehat{ r}_\kappa$. In
Figure~1, we have $\kappa = 6$. For an arbitrary fixed radius $r$
with $\widehat{ r}_\lambda < r < \widehat{ r}_{ \lambda + 1}$, and
some fixed $\lambda$ with $1 \leqslant \lambda \leqslant \kappa - 1$,
we consider all connected components $M_{ > r}^c$, $1 \leqslant c
\leqslant c_\lambda$, of the cut out hypersurface $M \cap \big\{ \vert
\! \vert z \vert \! \vert > r \big\}$. Their number $c_\lambda$ is the
same for all $r \in \big( \widehat{ r}_\lambda, \widehat{ r}_{ \lambda
+ 1} \big)$. In Figure~1, when $\widehat{ r}_ 3 < r < \widehat{ r}_4$,
we see three such components.

\smallskip
By descending discrete induction $r \mapsto r - \eta$, we show that
each such connected hypersurface $M_{ > r}^c \subset \big\{ \vert \!
\vert z \vert \! \vert > r \big\}$ bounds a certain domain
$\widetilde{ \Omega}_{ >r}^c \subset \big\{ \vert \! \vert z \vert \!
\vert > r \big\}$ which is relatively compact in $\C^n$ and that
holomorphic functions in $\mathcal{ V}_\delta ( M)$ do extend
holomorphically and uniquely to $\widetilde{ \Omega}_{ >r}^c $. While
approaching a lower Morse point, three different topological processes
will occur\footnote{A certain number of other simpler cases will also
happen, where the components $\widetilde{ \Omega}_{ >r}^c$ do grow
regularly with respect to holomorphic extension, possibly changing
topology.}: {\sl creating} a new component $\widetilde{ \Omega}_{ > r
- \eta}^{ c'}$ to be filled in further; {\sl merging} two components
$\widetilde{ \Omega}_{ > r - \eta}^{ c_1'}$ and $\widetilde{ \Omega}_{
> r - \eta}^{ c_2'}$ which meet; and {\sl suppressing} one superfluous
component $\widetilde{ \Omega}_{ > r - \eta }^{ c_1'}$.

\smallskip
The unavoidable multivaluation phenomenon will be tamed by the idea of
{\it separating ab initio}\, the components $M_{ > r}^c$, $1 \leqslant
c \leqslant c_\lambda$. Indeed, an advantageous topological property
will be shown to be inherited through the induction $r \mapsto r -
\eta$, hence always true, namely that two different domains
$\widetilde{ \Omega}_{ > r}^{ c_1}$ and $\widetilde{ \Omega}_{ > r}^{
c_1}$ are either disjoint or one is contained in the other.
Consequently, the multivaluation aspect will only happen in the sense
that the two {\it uniquely defined and univalent}\, holomorphic
extensions $f_{ r}^{ c_1}$ to $\widetilde{ \Omega}_{ > r}^{ c_1}$
and $f_{ r}^{ c_2}$ to $\widetilde{ \Omega}_{ > r}^{ c_2}$ can be
different on $\widetilde{ \Omega}_{ > r}^{ c_1}$, in case $\widetilde{
\Omega}_{ > r}^{ c_1} \subset \widetilde{ \Omega}_{ > r}^{ c_2}$, or
{\it vice versa}. In this way, we {\it avoid completely}\, to deal
with Riemann domains spread over $\C^n$.

\smallskip
Some of the elements of our approach should be viewed in a broader
context.  In their celebrated paper \cite{ag1962} ({\it see} also
\cite{hl1988}), Andreotti and Grauert observed that convenient
exhaustion functions can be used to prove very general extension and
finiteness results on $q$-concave complex varieties. Their arguments
implicitly contained a geometrical proof of the Hartogs extension
theorem in the case where the domain $\Omega \subset \C^n$ is
pseudoconvex (whence Forn{\ae}ss' counter-example {\it must}\, be
nonpseudoconvex).  However, in contrast to our finer method, the
existence of an internal strongly pseudoconvex exhaustion function
$\rho$ on a complex manifold $X$ excludes {\it ab initio}\,
multisheetedness: indeed, in such a circumstance, extension holds
stepwise from shells of the form $\big\{ z\in X: \, a < \rho (z) < b
\big\}$ just to deeper shells $\{ a' < \rho < b \}$ with $a ' < a$
(details are provided in~\cite{ mp2007}), namely the topology is
controlled in advance by $\rho$ and multiple domains as $\widetilde{
\Omega}_{ > r}^c$ above cannot at all appear.

There is a nice alternative approach to the (singular) Hartogs
extension theorem via a global continuity principle, realized in
\cite{jp2002} by J\"oricke and the second author, with the purpose of
understanding removable singularities by means of (geometric)
envelopes of holomorphy. The idea is to perform holomorphic extensions
along one-parameter families of holomorphic curves (not suppose to be
discs). A basic extension theorem on some appropriate Levi flat
3-manifolds, called {\em Hartogs manifolds}\, in \cite{jp2002}, is
shown via stepwise extension in the direction of an increasing real
parameter. The geometrical scheme of this construction has a common
topological element with our method: the simultaneous holomorphic
extension to collections of domains that are pairwise either disjoint
or one is contained in the other.

On the other hand, our technique only rely upon the existence of
appropriate exhaustion functions, without requiring neither the
existence of Levi-flat 3-manifolds nor the existence of global
holomorphic functions in the ambient complex manifold. In addition,
inspired by a definition formulated by Forn{\ae}ss in~\cite{ fo1998},
we establish that only a {\it finite}\, number of Levi-Hartogs figures
is needed in the filling process.  Finally, we would like to mention
that a straightforward adaptation of the proof developed here would
yield a geometrical proof of the Hartogs-type extension theorem of
Andreotti and Hill (\cite{ah1972}), which is valid for arbitrary
domains in $(n-1)$-complete manifolds (in the sense of
Andreotti-Grauert \cite{gr1994}).

\smallskip

Twenty-two colored illustrations appear, each one being inserted at
the appropriate place in the text. Abstract geometrical thought being
intrinsically pictural, we hope to address to a broad audience
of complex analysts and geometers.

\section*{ \S2.~Preparation of the boundary and unique extension}

\subsection*{ 2.1.~Preparation of a good $\mathcal{ C }^\infty$
boundary} Denote by $\vert \! \vert z \vert \! \vert := \big( \vert
z_1 \vert^2 + \cdots + \vert z_n \vert^2 \big)^{ 1/2}$ the Euclidean
norm of $z = (z_1, \dots, z_n) \in \C^n$ and by $\B^n ( p, \delta ) :=
\big\{ \vert \! \vert z - p \vert \! \vert < \delta \big\}$ the open
ball of radius $\delta >0$ centered at a point $p$. If $E \subset
\C^n$ is any set,
\[
\mathcal{V}_\delta(E)
:= 
\cup_{p\in E}\, 
\B^n(p,\delta)
\] 
is a concrete open neighborhood of $E$.

As in the Hartogs theorem, assume that the domain $\Omega \Subset
\C^n$ has connected boundary $\partial \Omega$ and let $\mathcal{ V}
(\partial \Omega)$ be an open neighborhood of $\partial \Omega$, also
connected. Clearly, there exists $\delta_1$ with $0 < \delta_1 < \! \!
< 1$ such that $\partial \Omega \subset \mathcal{ V}_{ \delta_1} (
\partial \Omega) \subset \mathcal{ V} (\partial \Omega)$; of course,
$\mathcal{ V}_{ \delta_1 } ( \partial \Omega)$ is then also
connected. Choose a point $p_0 \in \C^n$ with ${\rm dist} \, (p_0,
\overline{ \Omega }) = 3$, center the coordinates $(z_1,
\dots, z_n)$ at $p_0$ and consider the distance function 
\def\theequation{2.2}\begin{equation}
r(z) 
:=
\vert\!\vert
z-p_0 
\vert\!\vert
= 
\vert\!\vert
z 
\vert\!\vert. 
\end{equation}
It is crucial to prepare as follows the boundary, replacing $(\Omega,
\partial \Omega)$ by $(\Omega_M, M)$, thanks to some transversality
arguments that are standard in Morse theory {\rm (\cite{ mi1963} and
\cite{ hi1976}, Ch.~6)}.

\def\thelemma{2.3}\begin{lemma}
There exists a $\mathcal{ C}^\infty$ connected closed and oriented
hypersurface $M \subset \mathcal{ V}_{ \delta_1 /2} ( \partial
\Omega)$ such that{\rm :}

\begin{itemize}

\smallskip\item[{\bf (i)}]
$M$ bounds a unique bounded domain $\Omega_M$ with $\Omega \subset
\Omega_M \cup \mathcal{ V} (\partial \Omega)${\rm ;}

\smallskip\item[{\bf (ii)}]
the restriction $r_M (z) := r (z) 
\big\vert_M$ of the distance function $r (z) =
\vert \! \vert z \vert \! \vert$ to $M$ has only a finite
number $\kappa$ of critical points $\widehat{ p}_\lambda \in M$, $1
\leqslant \lambda \leqslant \kappa$, located on different sphere
levels, namely
\[
2\leqslant 
r(\widehat{p}_1)
<
\cdots
<
r(\widehat{p}_\kappa)
\leqslant
5
+
\text{\rm diam}(\overline{\Omega}){\rm ;}
\]

\smallskip\item[{\bf (iii)}]
all the $(2 n - 1) \times (2n - 1)$ Hessian matrices ${\sf H} [ r_M ]
( \widehat{ p }_1 ), \dots, {\sf H} [ r_M ] ( \widehat{ p }_\kappa )$
have a nonzero determinant.

\end{itemize}\smallskip
\end{lemma}

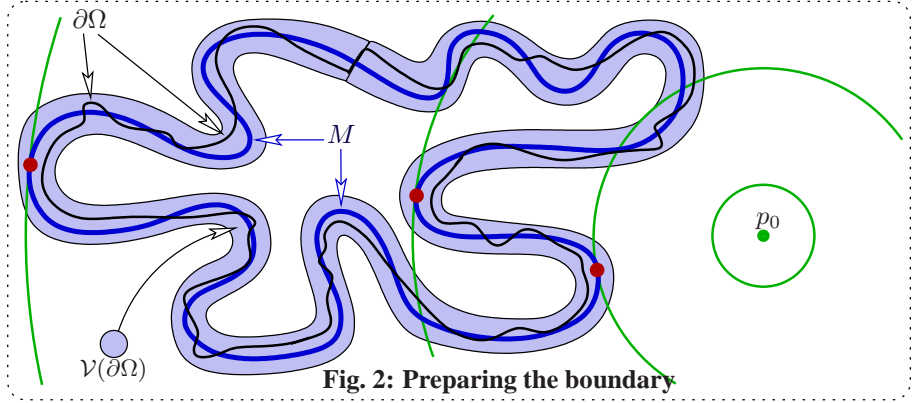
\begin{center}
\begin{picture}(0,0)%
\includegraphics{octopus.pstex}%
\end{picture}%
\setlength{\unitlength}{4144sp}%
\begingroup\makeatletter\ifx\SetFigFont\undefined
\def\x#1#2#3#4#5#6#7\relax{\def\x{#1#2#3#4#5#6}}%
\expandafter\x\fmtname xxxxxx\relax \def\y{splain}%
\ifx\x\y   
\gdef\SetFigFont#1#2#3{%
  \ifnum #1<17\tiny\else \ifnum #1<20\small\else
  \ifnum #1<24\normalsize\else \ifnum #1<29\large\else
  \ifnum #1<34\Large\else \ifnum #1<41\LARGE\else
     \huge\fi\fi\fi\fi\fi\fi
  \csname #3\endcsname}%
\else
\gdef\SetFigFont#1#2#3{\begingroup
  \count@#1\relax \ifnum 25<\count@\count@25\fi
  \def\x{\endgroup\@setsize\SetFigFont{#2pt}}%
  \expandafter\x
    \csname \romannumeral\the\count@ pt\expandafter\endcsname
    \csname @\romannumeral\the\count@ pt\endcsname
  \csname #3\endcsname}%
\fi
\fi\endgroup
\begin{picture}(5424,2409)(418,-1955)
\put(815,286){\makebox(0,0)[lb]{\smash{\SetFigFont{10}{12.0}{rm}$\partial\Omega$}}}
\put(2340,-409){\makebox(0,0)[lb]{\smash{\SetFigFont{10}{12.0}{rm}\blue{$M$}}}}
\put(4906,-893){\makebox(0,0)[lb]{\smash{\SetFigFont{10}{12.0}{rm}$p_0$}}}
\put(859,-1803){\makebox(0,0)[lb]{\smash{\SetFigFont{9}{10.8}{rm}$\mathcal{V}(\partial\Omega)$}}}
\put(2314,-1873){\makebox(0,0)[lb]{\smash{\SetFigFont{10}{12.0}{rm}{\bf Fig.~2: Preparing the boundary}}}}
\end{picture}

\end{center}

Sometimes, $r_M$ satisfying {\bf (ii)} and {\bf (iii)} is called a
{\sl good Morse function} on $M$. We will shortly say that $M$ is a
{\sl good boundary}.

If $k_\lambda$ is the number of positive eigenvalues of the
(symmetric) Hessian matrix ${\sf H} [ r_M ] ( \widehat{ p}_\lambda)$,
the extrinsic Morse lemma (\cite{ mi1963, hi1976}) shows that there
exist $2n$ real coordinates $\big( v, x_1, \dots, x_{ k_\lambda}, y_1,
\dots, y_{ 2n - k_\lambda - 1} \big)$ in a 
neighborhood of $\widehat{ p}_\lambda$ in $\C^n$ such that

\begin{itemize}

\smallskip\item[$\bullet$]
the sets $\{ v ( z) = {\rm cst} \}$ simply correspond\footnote{ In
fact, one can just take the translated radius $r(z) - r ( \widehat{
p}_\lambda)$ as the coordinate $v = v (z)$.} to the spheres $\{ r (z)
= {\rm cst} \}$ near $\widehat{ p}_\lambda$;

\smallskip\item[$\bullet$]
$\big( x_1, \dots, x_{ k_\lambda }, y_1, \dots, y_{ 2n - k_\lambda - 1}
\big)$ provide $( 2n -1)$ local coordinates on the hypersurface $M$,
whose graphed equation is normalized to be the simple hyperquadric
\[
v
=
\sum_{1\leqslant j\leqslant k_\lambda}\,
x_j^2
-
\sum_{1\leqslant j\leqslant 2n - k_\lambda-1}\,
y_j^2.
\]

\end{itemize}\smallskip

Classically, the number $(2n - k_\lambda - 1)$ of {\it negatives}\, is
called the {\sl Morse index} of $r(z) \big \vert_M$ at $\widehat{ p
}_\lambda$; we will call $k_\lambda$ its {\sl Morse coindex}.

For rather general differential-geometric objects, Morse theory
enables to control a significant part of homotopy groups and of
(co)homologies, {\it e.g.} via Morse inequalities. In our case, we
shall be able to control somehow the global topology of the cut-out
domains $\Omega_M \cap \{ \vert \! \vert z \vert \! \vert > r \}$ that
re external to closed balls of radius $r$, filling them progressively
by means of analytic discs contained in small (Levi-)Hartogs figures
(Section~3). We start by checking rigorously that the Hartogs theorem
can be reduced to some good boundary.

\subsection*{ 2.4.~Unique holomorphic extension} 
If $\mathcal{ U} \subset \C^n$ is open, $\mathcal{ O}
( \mathcal{ U })$ denotes the ring of holomorphic functions in
$\mathcal{ U}$.

\def\thedefinition{2.5}\begin{definition}{\rm Given two connected open
sets $\mathcal{ U}_1 \subset \C^n$ and $\mathcal{ U }_2 \subset \C^n$
with $\mathcal{ U}_1 \cap \mathcal{ U}_2$ nonempty, we will
say\footnote{ Because in the sequel, the union $\mathcal{ U}_1 \cup
\mathcal{ U}_2$ would sometimes be a rather long, complicated
expression ({\it see}~e.g.~\thetag{ 3.9}), hence uneasy to read, we
will also say that $\mathcal{ O} (\mathcal{ U }_1 )$ extends
holomorphically {\sf and uniquely} to $\mathcal{ U }_2$.} that {\sl
$\mathcal{ O} (\mathcal{ U }_1 )$ extends holomorphically to
$\mathcal{ U}_1 \cup \mathcal{ U }_2$} if {\rm :}
\begin{itemize}

\smallskip\item[$\bullet$]
the intersection $\mathcal{ U}_1 \cap \mathcal{ U}_2$ is
connected{\rm ;}

\smallskip\item[$\bullet$]
there exists an open nonempty set $\mathcal{ V} \subset \mathcal{ U}_1
\cap \mathcal{ U}_2$ such that for every $f_1 \in \mathcal{ O} (
\mathcal{ U}_1)$, there exist $f_2 \in \mathcal{ O} ( \mathcal{ U}_2)$
with $f_2 \vert_{ \mathcal{ V}} = f_1 \vert_{ \mathcal{ V }}$.
\end{itemize}\smallskip

}\end{definition}

It then follows from the principle of analytic continuation that $f_1
\vert_{ \mathcal{ U}_1 \cap \mathcal{ U}_2} = f_2 \vert_{ \mathcal{
U}_1 \cap \mathcal{ U}_2}$, so that the joint function $F$, equal to
$f_j$ on $\mathcal{ U}_j$ for $j=1, 2$, is well defined, is
holomorphic in $\mathcal{ U}_1 \cup \mathcal{ U}_2$ and extends $f_1$,
namely $F \vert_{ \mathcal{ U}_1} = f_1$.

In concrete extensional situations, the coincidence of $f_1$ with
$f_2$ is controlled only in some small $\mathcal{ V} \subset \mathcal{
U}_1 \cap \mathcal{ U}_2$, so the connectedness of $\mathcal{ U}_1
\cap \mathcal{ U}_2$ appears to be useful to insure monodromy.
Sometimes also, we shall briefly write $\mathcal{ O} ( \mathcal{ U}_1)
= \mathcal{ O} (\mathcal{ U}_1 \cup \mathcal{ U}_2) \big\vert_{
\mathcal{ U }_1}$, instead of spelling rigorously{\rm :}
\[
\forall\,f_1\in\mathcal{O}
\big(
\mathcal{U}_1
\big)
\ \ \ \ 
\exists\ F\in\mathcal{O}
\big(
\mathcal{U}_1\cup\mathcal{U}_2
\big)
\ \ \ \
\text{\rm such that}
\ \ \ \
F\big\vert_{\mathcal{U}_1}
=
f_1.
\]

\def\thelemma{2.6}\begin{lemma}
Suppose that for some $\delta$ with $0 < \delta \leqslant \delta_1 /2$
so small that $\mathcal{ V}_\delta (M) \simeq M \times ( -
\delta, \delta)$ is a thin tubular neighborhood of the good boundary
$M \subset \mathcal{ V}_{ \delta_1 /2} ( \partial \Omega) \subset
\mathcal{ V} (\partial \Omega)$, the Hartogs theorem holds for the
pair $( \Omega_M, \mathcal{ V}_\delta ( M ))${\rm :}
\[
\mathcal{O}
\big( 
\mathcal{V}_\delta( 
M)
\big)
=
\mathcal{O}
\big(
\Omega_M\cup
\mathcal{V}_\delta( 
M)
\big)
\big\vert_{\mathcal{V}_\delta(M)}.
\]
Then the general Hartogs extension property holds{\rm: }
\[
\mathcal{O}
\big( 
\mathcal{V}( 
\partial\Omega)
\big)
=
\mathcal{O}
\big(
\Omega\cup
\mathcal{V}( 
\partial\Omega)
\big)
\big\vert_{\mathcal{V}(\partial\Omega)}.
\]
\end{lemma}

\proof
Let $f \in \mathcal{ O} \big( \mathcal{ V} (\partial \Omega)
\big)$. By assumption, its restriction to $\mathcal{ V}_\delta (M)
\subset \mathcal{ V} ( \partial \Omega)$ enjoys an extension $F_\delta
\in \mathcal{ O} \big( \Omega_M \cup \mathcal{ V}_\delta (M)
\big)$. To ascertain that $f$ and $F_\delta$ coincide in $\Omega_M
\cap \mathcal{ V} ( \partial \Omega)$, connectedness 
of $\Omega_M
\cap \mathcal{ V} ( \partial \Omega)$
is welcome.

\begin{center}
\input connected-curve.pstex_t
\end{center}

Letting $p, q \in \Omega_M \cap \mathcal{ V} ( \partial \Omega)$,
there exists a $\mathcal{ C}^\infty$ curve $\gamma : [ 0, 1] \to
\mathcal{ V} ( \partial \Omega)$ connecting $p$ to $q$. If $\gamma$
meets $M$, let $p'$ be the first point on $\gamma \cap M$ and let $q'$
be the last one. We then modify $\gamma$, joining $p'$ to $q'$ by
means of a curve $\mu$ entirely contained in the connected
hypersurface $M$. It suffices to push $\mu$ slightly inside
$\Omega_M$ to get an appropriate curve running from $p$ to $q$ inside
$\Omega_M \cap \mathcal{ V} (\partial \Omega)$. Thus, $\Omega_M \cap
\mathcal{ V} ( \partial \Omega)$ is connected. It follows, moreover,
that the open set
\[
\big[
\Omega_M
\cup
\mathcal{V}_\delta(M)
\big]
\cap
\mathcal{V}(\partial\Omega)
=
\big[
\Omega_M
\cap
\mathcal{V}(\partial\Omega)
\big]
\bigcup
\mathcal{V}_\delta(M)
\]
is also connected, so the coincidence $f = F_\delta$, valid in
$\mathcal{ V}_\delta (M)$,
propagates to $\big[ \Omega_M \cap \mathcal{ V} (\partial \Omega)
\big] \cup \mathcal{ V}_\delta (M)$. Finally, the function
\[
F:=
\left\{
\aligned
F_\delta
\ \ \
&
\text{\rm in}
\ \ 
\Omega_M\cup\mathcal{V}_\delta(M),
\\
f
\ \ \
&
\text{\rm in}
\ \
\mathcal{V}(\partial\Omega)\backslash\overline{\Omega}_M,
\endaligned
\right.
\]
is well defined \big(since $F_\delta = f$ in $\mathcal{ V}_\delta (M)
\backslash \overline{ \Omega}_M \simeq M \times (0, \delta) $\big), is
holomorphic in
\[
\Omega_M\cup\mathcal{V}(\partial\Omega)
=
\Omega\cup\mathcal{V}(\partial\Omega)
\]
and coincides with $f$ in $\mathcal{ V} (\partial \Omega)$.
\endproof

Thus, we are reduced to establish global holomorphic extension
with some good, geometrically controlled data.

\def\thetheorem{2.7}\begin{theorem}
Let $M \Subset \C^n$ {\rm (}$n \geqslant 2${\rm )} be a {\rm
connected} $\mathcal{ C }^\infty$ hypersurface bounding a domain
$\Omega_M \Subset \C^n$. Suppose to fix ideas that $2 \leqslant {\rm
dist }\, \big( 0 , \overline{ \Omega }_M \big) \leqslant 5$ and assume
that the restriction $r_M := r\vert_M$ of the distance function $r (z)
= \vert \! \vert z \vert \! \vert$ to $M$
is a Morse function having only a finite number
$\kappa$ of critical points $\widehat{ p}_\lambda \in M$, $1 \leqslant
\lambda \leqslant \kappa$, located on different sphere levels{\rm :}
\[
2\leqslant 
\widehat{r}_1
:=
r(\widehat{p}_1)
<
\cdots
<
\widehat{r}_\kappa
:=
r(\widehat{p}_\kappa)
\leqslant
5
+
\text{\rm diam}
\big(\overline{\Omega}_M\big).
\]
Then there exists $\delta_1 >0$ such that for every $\delta$ with $0 <
\delta \leqslant 
\delta_1$, the {\rm (}tubular{\rm )} neighborhood $\mathcal{
V }_\delta (M)$ enjoys the global Hartogs extension property into 
$\Omega_M${\rm
:}
\[
\mathcal{O}
\big( 
\mathcal{V}_\delta( 
M)
\big)
=
\mathcal{O}
\big(
\Omega_M\cup
\mathcal{V}_\delta( 
M)
\big)
\big\vert_{\mathcal{V}_\delta(M)},
\]
by ``pushing'' analytic discs inside a finite number of Levi-Hartogs
figures {\rm (\S3.3)}, without using neither the
Martinelli kernel, nor solutions of an 
auxiliary $\overline{ \partial }$ problem.
\end{theorem}

\section*{ \S3.~Quantitative Hartogs-Levi extension
\\
by pushing analytic discs}

\subsection*{ 3.1.~The classical Hartogs figure}
Local Hartogs phenomena can now enter the scene. They involve
translating (``pushing'') analytic discs and they will provide
small, elementary extensional steps to fill in $\Omega_M$.

Given $\varepsilon \in \R$ with $0 < \varepsilon < \! \! < 1$ and $a
\in \N$ with $1\leqslant a \leqslant n - 1$, we split the coordinates
$z\in \C^n$
as $(z_1, \dots, z_a)$ together with $( z_{ a+1}, \dots, z_n)$, and we
define the {\sl $(n - a)$-concave Hartogs figure} by
\[
\aligned
\mathcal{H}_\varepsilon^{n-a}
&
:=
\Big\{
\max_{1\leqslant i\leqslant a}\,
\vert z_i\vert<1,\ \
\max_{a+1\leqslant j\leqslant n}\,
\vert z_j\vert
<
\varepsilon
\Big\}
\\
&
\ \ \ \ \ \ \ \
\bigcup
\Big\{
1-\varepsilon
<
\max_{1\leqslant i\leqslant a}\,
\vert z_i\vert<1,\ \
\max_{a+1\leqslant j\leqslant n}\,
\vert z_j\vert
<
1
\Big\}.
\endaligned
\]

\begin{center}
\input hartogs-figure.pstex_t
\end{center}

\def\thelemma{3.2}\begin{lemma}
$\mathcal{ O} \big( \mathcal{ H }_\varepsilon^{ n-a} \big)$ extends
holomorphically to the unit polydisc
\[
\widehat{\mathcal{H}}_\varepsilon^{n-a}
:= 
\big\{
z\in\C^n:\, 
\max_{1\leqslant i\leqslant n}\, 
\vert 
z_i 
\vert 
< 
1 
\big\}
=
\Delta^n.
\]
\end{lemma}

\proof
As in the diagram, we consider only $n=2$, $a=1$, the general case
being similar. Pick an arbitrary $f \in \mathcal{ O} \big( \mathcal{
H}_\varepsilon^{ 2-1 } \big)$. Letting $\varepsilon '$ with $0 <
\varepsilon' < \varepsilon$, letting $z_2 \in \C$ with $\vert z_2
\vert < 1$, the analytic disc
\[
\zeta
\longmapsto
\big(
[1-\varepsilon']\,\zeta,z_2
\big)
=:
A_{z_2}^{\varepsilon'}(\zeta),
\]
where $\zeta$ belongs to the closed 
unit disc $\overline{ \Delta} = \{ \vert \zeta \vert
\leqslant 1 \}$, has its boundary $A_{ z_2}^{ \varepsilon'} (\partial \Delta)
= A_{ z_2}^{ \varepsilon' } \big(\{ \vert \zeta \vert = 1 \} \big)$
contained in $\mathcal{ H }_\varepsilon^{ 2-1 }$, the set where $f$ is
defined. Lowering dimensions by a unit, we draw discs as (green)
segments and boundaries of discs as (green) bold points. Thus, we may
compute the Cauchy integral
\[
F(z_1,z_2)
:=
\frac{1}{2\pi i}\,
\int_{\partial\Delta}\,
\frac{f
\big(
A_{z_2}^{\varepsilon'}(\zeta)
\big)}{\zeta-z_1}\,d\zeta.
\]
Differentiating under the sum, the function $F$ is seen to be
holomorphic. In addition, for $\vert z_2 \vert < \varepsilon$, it
coincides with $f$, because the full closed disc $A_{ z_2}^{
\varepsilon '} \big( \overline{ \Delta} \big)$ is contained in
$\mathcal{ H}_\varepsilon^{ 2 - 1}$ and thanks
to Cauchy's formula. Clearly, the
$A_{ z_2 }^{ \varepsilon'} ( \Delta )$ all together fill in the bidisc
$\Delta^2$. One may think that, as $z_2$ varies, discs are
``pushed'' gently by a virtual thumb.
\endproof

\subsection*{ 3.3.~Levi extension and the Levi-Hartogs figure}
Geometrically, the standard Hartogs figure is not best suited to
perform holomorphic extension from a strongly (pseudo)concave boundary.
For instance, in the proof of Theorem~2.7, we will encounter
complements in $\C^n$ of some closed balls whose radius decreases step
by step, and more generally spherical shells whose thickness increases
interiorly. Thus, we delineate an appropriate set up.

For $r\in \R$ with $r > 1$ and for $\delta \in \R$ with $0 < \delta <
\! \! < 1$, the sphere ${\sf S}_r^{ 2n - 1} = \{ z\in \C^n : \, \vert
\! \vert z \vert \! \vert = r \}$ of radius $r$ is the interior (and
strongly concave) boundary component of the spherical shell domain
\[
\mathcal{S}_r^{r+\delta}
:=
\big\{
r<\vert\!\vert z\vert\!\vert<r+\delta
\big\}
=
\bigcup_{p\in {\sf S}_r^{2n-1}}\,
\B^n(p,\delta)
\cap
\{\vert\!\vert z\vert\!\vert>r\}.
\]

\begin{center}
\input levi-hartogs.pstex_t
\end{center}

Near a point $p \in {\sf S}_r^{ 2n - 1}$ (left figure), all copies of
$\C^{ n-1 }$ (in green) which are parallel to the complex tangent
plane $T_p^c {\sf S}_r^{ 2n -1 }$ and which lie above the real plane
$T_p {\sf S}_r^{ 2n -1 }$ are entirely contained in $\C^n \big
\backslash \overline{ \B }_r^n$. To remain inside the shell $\mathcal{
S}_r^{ r + \delta}$, we could (for instance) restraint our considerations to
some half-cylinder of diameter $\approx \delta$, but it will be better
to shape a convenient half parallelepiped. Accordingly, for two small
$\varepsilon_j >0$, $j = 1, 2$, we introduce a geometrically relevant
{\sl Levi-Hartogs figure} (right illustration, reverse
orientation):
\[
\aligned
\mathcal{LH}_{\varepsilon_1,\varepsilon_2}
&
:=
\Big\{
\max_{1\leqslant i\leqslant n-1}\,
\vert z_i\vert
<
\varepsilon_1,\ \ \
\vert 
x_n
\vert
<
\varepsilon_1,
\ \ \ 
-\varepsilon_2
<
y_n
<
0
\Big\}
\\
&
\ \ \ \ \ \ \ 
\bigcup\
\Big\{
\varepsilon_1
-
(\varepsilon_1)^2
<
\max_{1\leqslant i\leqslant n-1}\,
\vert z_i\vert
<
\varepsilon_1,
\ \ \
\vert
x_n
\vert
<
\varepsilon_1
\ \ \
\vert
y_n
\vert
<
\varepsilon_2
\Big\}.
\endaligned
\]
To fill in this (bed-like) figure, we just compute the Cauchy integral
on appropriate analytic discs (the (green) horizontal ones) whose
boundaries remain in $\mathcal{ LH}_{ \varepsilon_1, \varepsilon_2 }$.

\def\thelemma{3.4}\begin{lemma}
$\mathcal{ O} \big( \mathcal{ LH }_{ \varepsilon_1, \varepsilon_2 }
\big)$ extends holomorphically to the full parallelepiped
\[
\widehat{\mathcal{LH}_{\varepsilon_1,\varepsilon_2}}
:=
\Big\{
\max_{1\leqslant i\leqslant n-1}\,
\vert
z_i
\vert
<
\varepsilon_1,
\ \ \ 
\vert
x_n
\vert
<
\varepsilon_1,
\ \ \ 
\vert
y_n
\vert
<
\varepsilon_2
\Big\}.
\]
\end{lemma}

Next, we must reorient and scale $\mathcal{ LH }_{ \varepsilon_1,
\varepsilon_2 }$ in order to put it inside the shell. For every point
$p \in {\sf S}_r^{ 2n -1}$, there exists some complex unitarian affine
map
\[
\Phi_p:
\ \ \ 
z 
\longmapsto 
p 
+
Uz,
\]
with $U \in {\sf SU } (n, \C)$, sending the origin $0 \in \overline{
\mathcal{ LH } }_{ \varepsilon_1, \varepsilon_2 }$ to $p$ and $T_0
\mathcal{ LH }_{ \varepsilon_1, \varepsilon_2 }$ to $T_p {\sf S}_r^{
2n -1}$, which in addition sends the half-parallelepiped (open) part
outside $\overline{ \B}_r^n$. But we have to insure that $\Phi_p \big(
\mathcal{ LH}_{ \varepsilon_1, \varepsilon_2} \big)$ {\it as a whole}
(including the thin walls) lies outside $\overline{ \B }_r^n$.

\def\thelemma{3.5}\begin{lemma}
If $\varepsilon_1 = c\, \delta$ and $\varepsilon_2 = c\, \delta^2$
with some appropriate\footnote{ We let the letter $c$ (resp. $C$)
denote a positive constant $<1$ (resp. $> 1$), absolute or depending
only on $n$, which is allowed to vary with the context.} positive
constant $c < 1$, then $\Phi_p \big( \mathcal{ LH}_{ \varepsilon_1,
\varepsilon_2 } \big)$ is entirely contained in the shell $\mathcal{
S}_r^{ r + \delta }$. Furthermore, $\Phi_p \big( \widehat{ \mathcal{
LH }_{ \varepsilon_1, \varepsilon_2 }} \big)$ contains a rind of
thickness $c\, \frac{ \delta^2 }{ r}$ around some region ${\sf R}_p
\subset {\sf S}_r^{ 2n -1 }$ whose $(2n -1)$-dimensional area equals
$\simeq c \, \delta^{ 2n - 1}$.
\end{lemma}

\smallskip

\begin{center}
\input green-lemon-rind.pstex_t
\end{center}

By a (radial) {\sl rind of thickness $\eta >0$ around an 
open region ${\sf R} \subset {\sf S}_r^{ 2n -1}$}, we mean
\[
{\sf Rind}
\big(
{\sf R},\eta
\big)
:=
\big\{
(1+s)z:\
z\in{\sf R},\
\vert s\vert
<
\eta/r
\big\}.
\]
We require that $\vert s \vert < \eta / r$ to insure that at every $z
\in {\sf R}$, the half-line $(0 z)^+$ emanating from the origin
intersects ${\sf Rind} \big( {\sf R}, \eta \big)$ along a 
symmetric segment of
length $2 \, \eta$ centered at $z$.

In the diagram above, we draw (in green) only the lower part of the
small region ${\sf R}_p$ got in Lemma~3.5. Its shape, when projected
onto $T_p {\sf S}_r^{ 2n -1}$, can either be (approximately) a
parallelepiped $\big\{ \vert z' \vert < c\, \delta, \ \vert x_n \vert
< c\, \delta \big\}$, as in the figure, or say, a ball $\big\{ \big(
\vert \! \vert z' \vert \! \vert^2 + \vert x_n \vert^2 \big)^{ 1/2} <
c\, \delta \big\}$; only the scaling constant $c$ changes.

The rigorous proof of the lemma (not developed here) involves
elementary reasonings with geometric inequalities and a dry explicit
control of the constants that does not matter for the sequel. 
The main argument uses the fact that ${\sf S }_r^{ 2n - 1}$ detaches
quadratically from $T_p {\sf S}_r^{ 2n -1}$, similarly as the parabola
$\big\{ y = - \frac{ 1}{ r }\, x^2 \big\}$ separates from the line $\{
y = 0 \}$ in $\R_{ x, y}^2$.

\smallskip

Since the area of ${\sf S}_r^{ 2n -1}$ equals $\frac{ 2\, \pi^n}{
(n-1) ! }\, r^{ 2 n -1} = C \, r^{ 2 n -1}$, by covering ${\sf S}_r^{
2n -1}$ with such adjusted ${\sf R}_p \subset \Phi_p \big( \widehat{
\mathcal{ LH}_{ \varepsilon_1, \varepsilon_2 }} \big)$ of area $c\,
\delta^{ 2n -1}$ and by controlling monodromy ({\it see} rigorous
arguments below) we deduce:

\def\thecorollary{3.6}\begin{corollary}
By means of a finite number $\leqslant C \, \big( \frac{ r}{ \delta}
\big)^{ 2 n - 1 }$ of Levi-Hartogs figures, $\mathcal{ O} \big(
\mathcal{ S }_r^{ r + \delta } \big)$ extends holomorphically to the
slightly deeper spherical shell $\mathcal{ S}_{ r - c\, \frac{
\delta^2 }{ r}}^{ r + \delta }$.
\end{corollary}

This application could seem superfluous, because large analytic discs
with boundaries contained in $\mathcal{ S }_r^{ r + \delta }$ would
yield holomorphic extension to the whole ball $\B_{ r + \delta }^n$ in
one single step. However, in our situation illustrated by Figure~1,
when intersecting ${\sf S }_r^{ 2n -1}$ with the neighborhood
$\mathcal{ V }_\delta ( M )$, we shall only get small subregions of
${\sf S}_r^{ 2n -1 }$. Hopefully, thanks to our local Levi-Hartogs
figures, we may obtain a suitable semi-global extensional statement,
valuable for proper subsets of the shell $\mathcal{ S}_r^{ r + \delta
}$ whose shape is arbitrary. The next statement, not available by
means of large discs, will be used a great number of times
in the sequel.

\def\theproposition{3.7}\begin{proposition}
Let ${\sf R} \subset {\sf S}_r^{ 2n -1}$ {\rm (}with $r >1$ and $n
\geqslant 2${\rm )} be a relatively open set having $\mathcal{
C}^\infty$ boundary ${\sf N} := \partial {\sf R}$ and let $\delta >0$
with $0 < \delta < \! \! < 1$. Then holomorphic functions in the open
piece of shell {\rm (}a one-sided neighborhood of ${\sf R} \cup {\sf
N}${\rm )}{\rm :}
\[
\aligned
\text{\rm Shell}_r^{r+\delta}
\big(
{\sf R}\cup{\sf N}
\big)
:=
&\
\big(
\C^n
\big\backslash
\overline{\B}_r^n
\big)
\cap
\mathcal{V}_\delta
\big(
{\sf R}\cup{\sf N}
\big)
\\
=
&\
\bigcup_{p\in{\sf R}\cup{\sf N}}\,
\B^n(p,\delta)
\cap
\{
\vert\!\vert
z
\vert\!\vert
>
r
\}
\endaligned
\]
do extend holomorphically to a rind of thickness $c\, \frac{ \delta^2
}{ r}$ around ${\sf R}$ by means of a finite number $\leqslant C \,
\frac{ \text{ \rm area} ( {\sf R }) }{ \delta^{ 2n-1 }}$ of
Levi-Hartogs figures.
\end{proposition}

\begin{center}
\input buoy.pstex_t
\end{center}

\proof
We must control uniqueness of holomorphic extension (monodromy) into
rinds covered by successively attached Levi-Hartogs figures.
Noticing $c\,\delta^2\,r^{ -1} < \! \! < \delta$, 
the considered rinds are much thinner than the
piece of shell.

\def\thelemma{3.8}\begin{lemma}
If ${\sf R}' \subset {\sf R}$ is an arbitrary open subset and if ${\sf
R}_{ p'} \subset \Phi_{p'} 
\big( \widehat{ \mathcal{ LH}_{ \varepsilon_1,
\varepsilon_2}} \big)$ is a small Levi-Hartogs region centered at an
arbitrary point $p' \in {\sf R}$, then the intersection
\def\theequation{3.9}\begin{equation}
{\sf Rind}
\big(
{\sf R}_{p'},\,c\,\delta^2\,r^{-1}
\big)
\bigcap
\Big(
\text{\rm Shell}_r^{r+\delta}
\big(
{\sf R}\cup{\sf N}
\big)
\bigcup
{\sf Rind}
\big(
{\sf R}',\,c\,\delta^2\,r^{-1}
\big)
\Big)
\end{equation}
is connected.
\end{lemma}

Admitting the lemma for a while, we pick a finite number $m \leqslant
C \, \frac{ \text{\rm area} ( {\sf R}) }{ \delta^{ 2n -1}}$ of points
$p_1, \dots, p_m \in {\sf R} \cup {\sf N}$ such that the associated
local regions ${\sf R}_{ p_k}$ contained in the filled Levi-Hartogs
figures $\Phi_{ p_k} \big(
\widehat{ \mathcal{ LH}_{ \varepsilon_1, 
\varepsilon_2}} \big)$
provided by Lemma~3.5 do cover ${\sf R} \cup {\sf N}$, namely
${\sf R}_{ p_1} \cup \cdots \cup 
{\sf R}_{ p_m} \supset {\sf R} \cup {\sf N}$.

Starting with ${\sf R}' := \emptyset$ and $p' := p_1$, unique
holomorphic extension of $\mathcal{ O} \big( \text{\rm Shell}_r^{ r+
\delta} ( {\sf R} \cup {\sf N} ) \big)$ to ${\sf Rind} \big( {\sf
R}_{ p'}, c\, \delta^2 \, r^{ -1} \big)$ holds by means of Lemma~3.4,
monodromy being assured thanks to the connectedness of the
intersection~\thetag{ 3.9}. Reasoning
by induction, fixing some $k$ with
$1\leqslant k \leqslant m -1$, setting ${\sf R}' := \cup_{ 1\leqslant
j \leqslant k}\, {\sf R}_{ p_j}$, $p' := p_{ k+1}$ and assuming that
unique holomorphic extension is got from $\text{\rm Shell}_r^{ r+
\delta} \big( {\sf R} \cup {\sf N} \big)$ into
\begin{small}
\[
\text{\rm Shell}_r^{r+\delta}
\big({\sf R}\cup{\sf N}\big)
\bigcup
{\sf Rind}
\big(
{\sf R}',c\,\delta^2\,r^{-1}
\big)
=
\text{\rm Shell}_r^{r+\delta}
\big({\sf R}\cup{\sf N}\big)
\bigcup_{1\leqslant j\leqslant k}\,
{\sf Rind}
\big(
{\sf R}_{p_j},c\,\delta^2\,r^{-1}
\big),
\]
\end{small}

\noindent
we add the Levi-Hartogs figure $\Phi_{ p_{ k+1}} \big( \widehat{
\mathcal{ LH}_{ \varepsilon_1, \varepsilon_2}} \big)$ constructed in
Lemma~3.5, and we get unique holomorphic extension to ${\sf Rind}
\big( {\sf R}_{ p_{ k+1}}, c\, \delta^2\, r^{ -1}\big)$, monodromy
being assured again thanks to the connectedness of the
intersection~\thetag{ 3.9}. Since ${\sf Rind} \big( {\sf R}, c\,
\delta^2 \, r^{ -1} \big) \subset \bigcup_{1\leqslant k \leqslant m}\,
{\sf Rind} \big( {\sf R}_{ p_k}, c\, \delta^2 \, r^{ -1} \big)$, the
proposition is proved.
\endproof

\proof[Proof of Lemma~3.8]
To establish connectedness of the open set~\thetag{ 3.9}, picking two
arbitrary points $q_0$, $q_1$ in it, we must produce a curve joining
$q_0$ to $q_1$ inside~\thetag{ 3.9}. The two radial segments of
length $2\, c\, \delta^2 \, r^{ -1}$ passing through $q_0$ and $q_1$
that are 
centered at two appropriate points of ${\sf S}_r^{ 2n -1}$ are by
definition both entirely contained in ${\sf Rind} \big( {\sf R}_{
p'}, c\, \delta^2 \, r^{ -1} \big)$ as well as in ${\sf Rind}
\big( {\sf R}', c\, \delta^2 \, r^{ -1} \big)$. Thus, moving radially,
we may join inside~\thetag{ 3.9} $q_0$ to a new point $q_0'$ and $q_1$
to a new point $q_1'$, which both belong to the upper half-rind
\[
\big\{
(1+s)\,z:\
z\in{\sf R}_{p'},\
0<s<c\,\delta^2\,r^{-1}\big/r
\big\}.
\]
Since this upper half-rind is connected and contained in $\text{\rm
Shell}_r^{ r+ \delta} \big( {\sf R} \cup {\sf N} \big)$, we may
finally join inside~\thetag{ 3.9} the point $q_0'$ to $q_1'$.
\endproof

In the sequel, in order to avoids several gaps and traps, we will put
emphasis on rigourously checking univalence of holomorphic extensions.

\section*{ \S4.~Filling domains outside balls of decreasing radius}

\subsection*{ 4.1.~Global Levi-Hartogs filling from the farthest
point} We can now launch the proof of Theorem~2.7. The $\delta_1$ is
first chosen so small that $\mathcal{ V}_\delta ( M)$ is a true
tubular neighborhood of $M$ for every $\delta$ with $0 < \delta
\leqslant \delta_1$. Shrinking even more $\delta_1$, in balls of
radius $\delta_1$ centered at its points, the hypersurface $M$ is well
approximated by its tangent planes.

The farthest point of $\overline{ \Omega }_M$ from the origin is unique
and it coincides with $\widehat{ p }_\kappa$ since by assumption
$\widehat{ p }_\kappa$ is the single critical point of $r ( z)
\big \vert_M$ with $\vert \! \vert \widehat{ p}_\kappa \vert \! \vert =
\max_{ 1\leqslant \lambda \leqslant \kappa} \, \vert \! \vert
\widehat{ p}_\lambda \vert \! \vert$. By assumption also, the
Hessian matrix of $r ( z) \big \vert_M$ is nondegenerate at $\widehat{
p }_\kappa$; this also follows automatically from the inclusion
$\overline{ \Omega }_M \subset \overline{ \B }_{ \widehat{ r
}_\kappa}^n$, which constrains strong convexity of $M$ at $\widehat{ p
}_\kappa$. Consequently, according to the {\sl Morse lemma} (\cite{
mi1963}, \cite{ hi1976}, Ch.~6), 
there exist local coordinates $( \theta_1, \dots,
\theta_{ 2n - 1 })$ on $M$ centered at $\widehat{ p }_\kappa$ such
that the intersection $M \cap {\sf S }_r^{ 2n - 1 }$ is given by the
equation
\[
-\theta_1^2
-\cdots
-\theta_{2n-1}^2
=
r
-
\widehat{r}_\kappa,
\]
for all $r$ close to $\widehat{ r }_\kappa$. Thus $M \cap {\sf S}_r^{ 2n
-1}$ is empty for $r > \widehat{ r }_\kappa$; it reduces to $\{
\widehat{ p }_\kappa \}$ for $r = \widehat{ r }_\kappa$; and it is
diffeomorphic to a $( 2n - 2)$-sphere for $r < \widehat{ r }_\kappa$
close to $\widehat{ r}_\kappa$.

Similarly, the nearest point of $\overline{ \Omega }_M$ from the
origin is unique and it coincides with $\widehat{ p}_1$; notice that
hence $\kappa \geqslant 2$. Also, the second farthest critical point
$\widehat{ p}_{ \kappa - 1}$ lies at a distance $\widehat{ r}_{ \kappa
- 1} < \widehat{ r }_\kappa$ from $0$. If necessary, we shrink
$\delta_1$ to insure
\def\theequation{4.2}\begin{equation}
\delta_1 
< 
\!\!
< 
\min_{1 \leqslant \lambda \leqslant \kappa - 1} 
\big\{ 
\widehat{r}_{\lambda+1} 
- 
\widehat{r}_\lambda 
\big\}.
\end{equation}
Next, for every radius $r$ with $\widehat{ r}_{ \kappa - 1} < r <
\widehat{ r}_\kappa$, we introduce the cut out domain
\[
\Omega_{>r}
:=
\Omega_M
\cap
\big\{
\vert\!\vert
z
\vert\!\vert
>
r
\big\}
\]
together with the cut out hypersurface
\[
M_{>r}
:=
M
\cap
\big\{
\vert\!\vert
z
\vert\!\vert
>
r
\big\}.
\]

\begin{center}
\input farthest-point.pstex_t
\end{center}

Since there are no critical points of $r( z) \big\vert_M$ in the
interval $\big( \widehat{ r}_{ \kappa - 1}, \widehat{ r}_\kappa
\big)$, Morse theory shows that $M_{ >r}$ is a deformed spherical cap
diffeomorphic to $\R^{ 2n -1}$ for every $r$ with $\widehat{ r}_{
\kappa -1} < r < \widehat{ r}_\kappa$. Also, $\Omega_{ >r}$ is then a
piece of deformed ball diffeomorphic to $\R^{ 2 n }$.

The boundary in $\C^n$ of $\Omega_{ >r }$
\[
\partial\Omega_{>r}
=
M_{>r}
\cup
{\sf R}_r
\cup
{\sf N}_r
\]
consists of $M_{ >r}$ together with the open subregion ${\sf R}_r :=
\Omega_M \cap \big\{ \vert \! \vert z \vert \! \vert = r \big\}$ of
${\sf S}_r^{ 2n -1}$ which is diffeomorphic to $\R^{ 2n -1}$ and has
boundary ${\sf N}_r := M \cap \big\{ \vert \! \vert z \vert \! \vert =
r \big\}$ diffeomorphic to the unit $( 2n - 2)$-sphere. Thus, the
global geometry of $\Omega_{ >r}$ is understood.

We can also cut out $\mathcal{ V}_\delta (M)$, getting $\mathcal{
V}_\delta (M)_{ > r}$. The central figure shows that when $r >
\widehat{ r}_{ \kappa - 1}$ is very close to $\widehat{ r}_{ \kappa
-1}$, a parasitic connected component $\mathcal{ W}_{ >r}$ of
$\mathcal{ V}_\delta ( M)_{ >r}$ might appear near $\widehat{ p}_{
\kappa -1}$. After filling $\Omega_{> r}$ progressively by means of
Levi-Hartogs figures ({\it see} below), because $\Omega_{ >r} \cap
\mathcal{ V}_\delta ( M)_{ > r}$ is {\it not}\, connected in such a
situation, {\it no}\, unique holomorphic extension can be assured, and
in fact, multivalence might well occur.

A trick to erase such parasitic components $\mathcal{ W}_{ >r}$ is to
consider instead the open set 
\[
\mathcal{V}_\delta
\big(M_{>r}\big)_{>r}
=
\mathcal{V}_\delta
\big(M_{>r}\big)\cap
\big\{ 
\vert\!\vert 
z 
\vert\!\vert
> 
r 
\big\},
\] 
putting a double ``${}_{>r}$''. It is drawn in the right
figure and it is always diffeomorphic to $M_{ >r} \times (- \delta,
\delta)$.

From pieces of shells as in Proposition~3.7 which embrace spheres of
varying radius $r$, holomorphic extension holds to (symmetric) rinds
whose thickness $c\, \delta\, r^{ -1}$ also varies. To simplify, we
introduce the smallest appearing thickness
\def\theequation{4.3}\begin{equation}
\eta
:=
\min_{\widehat{r}_1\leqslant r\leqslant\widehat{r}_\kappa}\,
c\,\delta\,r^{-1}
=
c\,\delta\,\widehat{r_\kappa}^{-1}, 
\end{equation}
and we observe that it follows trivially from Proposition~3.7 (just by
shrinking and by restricting) that holomorphic extension holds to some
rind around ${\sf R}$ of {\it arbitrary}\, smaller thickness $\eta'
>0$ with $0 < \eta' \leqslant \eta$. In the sequel, our rinds shall
most often have the uniform thickness $\eta$, and sometimes also, a
smaller one $\eta'$. Shrinking the constant $c$ of $\eta$
in~\thetag{ 4.3}, we
insure $\eta < \! \! < \delta_1$.

Summarizing, we list and we compare the quantities introduced so far:
\begin{small}
\def\theequation{4.4}\begin{equation}
\left\{
\aligned
&
0
<
\delta
\leqslant
\delta_1
\ \ \ \ \ \ \ \ \ \ \ \ \ \ \ \ \ \
\ \ \ \ \ \ \ \ \ \ \ \ \ \ \ \ \ \
\ \ \ \ \ \ \ \ \ \ \ \ \ \ \ \ \ \
\text{\rm neighborhood}
\ 
\mathcal{V}_\delta(M) 
\\
&
2\leqslant 
r(\widehat{p}_1)
<
\cdots
<
r(\widehat{p}_\kappa)
\leqslant
5
+
\text{\rm diam}
\big(\overline{\Omega}_M\big)
\ \ \ \ \ \ \ \ \ \ \ \
\text{\rm Morse radii}
\\
&
\delta
\leqslant
\delta_1 
< 
\!\!
< 
\min_{1 \leqslant \lambda \leqslant \kappa - 1} 
\big\{ 
\widehat{r}_{\lambda+1} 
- 
\widehat{r}_\lambda 
\big\}
\ \ \ \ \ \ \ \ \ \ \ \ \ \ \ \ \ \ \ \ \
\text{\rm smallness of}
\ 
\mathcal{V}_\delta(M) 
\\
&
\eta
:=
c\,\delta^2\,\widehat{r}_\kappa^{-1}
\ \ \ \ \ \ \ \ \ \ \ \ \ \ \ \ \ \ \
\ \ \ \ \ \ \ \ \ \ \ \ \ \ \ \ \ \ \ \ \
\text{\rm uniform useful rind thickness}
\\
&
\eta
<\!\!<
\delta
\ \ \ \ \ \ \ \ \ \ \ \ \ \ \ \ \ \ \
\ \ \ \ \ \ \ \ \ \ \ \ \ \ \ \ \ \ \ \ 
\text{\rm thickness of extensional rinds is tiny}
\endaligned\right.
\end{equation}
\end{small}

\def\theproposition{4.5}\begin{proposition}
For every cutting radius $r$ with $\widehat{ r }_{ \kappa -1 } < r <
\widehat{ r }_\kappa$ arbitrarily close to $\widehat{ r }_{ \kappa -1 }$,
holomorphic functions in the open set
\[
\mathcal{V}_\delta\big( M_{>r} \big)_{>r}
=
\mathcal{V}_\delta
\big(M_{>r}\big)\cap
\big\{ 
\vert\!\vert 
z 
\vert\!\vert
> 
r 
\big\}
\]
do extend holomorphically and uniquely to $\Omega_{ >r }$ by means of
a finite number $\leqslant C\, \big( \frac{ \widehat{ r}_\kappa}{
\delta} \big)^{ 2n -1} \big[ \frac{ \widehat{ r}_\kappa - r }{ \eta}
\big]$ of Levi-Hartogs figures.
\end{proposition}

\proof
We fix such a radius $r$ with $\widehat{ r}_{ \kappa -1} < r <
\widehat{ r}_\kappa$. Putting a single Levi-Hartogs figure at
$\widehat{ p }_\kappa$ as in Proposition~3.7, we get unique
holomorphic extension to $\Omega_{ > \widehat{ r }_\kappa - \eta
}$. Since $\eta < \! \! < \delta$, we have $\widehat{ r }_\kappa -
\eta > \widehat{ r}_{ \kappa -1 }$. If the radius $\widehat{ r
}_\kappa - \eta$ is already $< r$, we just shrink to $\eta ':=
\widehat{ r}_\kappa - r < \eta$ the thickness of our single
rind, getting unique holomorphic extension to $\Omega_{ > \widehat{
r}_\kappa - \eta '} = \Omega_{ > r}$.

Performing 
induction on an auxiliary integer $k \geqslant 1$, we suppose that,
by descending from $\widehat{ r}_\kappa$ to a lower radius $r ' :=
\widehat{ r }_\kappa - k \eta$ assumed to be still $\geqslant 
r$, holomorphic
functions in $\mathcal{ V }_\delta \big( M_{ >r } \big)_{ >r}$ extend
holomorphically {\it and uniquely} (remind Definition~2.5) to $\Omega_{ >
r'}$.

\def\thelemma{4.6}\begin{lemma}
For every radius $r'$ with $\widehat{ r}_{ \kappa -1} < r < r' <
\widehat{ r }_\kappa$,
\def\theequation{4.7}\begin{equation}
\text{\rm Shell}_{r'}^{r'+\delta}
\big(
{\sf R}_{r'}
\cup 
{\sf N}_{r'}
\big)
\ \
\text{\rm is contained in}
\ \
\Omega_{>r'}
\bigcup
\mathcal{V}_{\delta}
\big(
M_{>r}
\big)_{>r}.
\end{equation}
\end{lemma}

\begin{center}
\input cap-shell.pstex_t
\end{center}

\proof
Picking an arbitrary point $p \in {\sf R}_{ r'} \cup {\sf N}_{ r'}$, we
must verify that 
\[
\B^n(p,\delta) 
\cap
\{\vert\!\vert 
z 
\vert\!\vert 
> 
r'
\}
\] 
is contained in the right hand side of~\thetag{ 4.7}.

If $p \in {\sf N}_{ r'} \subset 
M$, whence $p \in M_{ >r}$, we get
simply what we want{\rm :}
\[
\aligned
\B^n(p,\delta)
\cap
\{
\vert\!\vert z\vert\!\vert>r'
\}
&
\subset
\mathcal{V}_\delta
\big(
M_{>r}
\big)
\cap
\{
\vert\!\vert z\vert\!\vert>r'
\}
\\
&
\subset
\mathcal{V}_\delta
\big(
M_{>r}
\big)
\cap
\{
\vert\!\vert z\vert\!\vert>r
\}
\\
&
=
\mathcal{V}_\delta
\big(
M_{>r}
\big)_{>r}.
\endaligned
\]

If $p \in {\sf R}_{ r'} \big \backslash {\sf N}_{ r'}$, whence $p
\in \Omega_M$, reasoning by contradiction, we assume that there exists
a point $q \in \B^n ( p, \delta) \cap \{ \vert \! \vert z \vert \!
\vert > r' \}$ in the cut out ball which does not belong to the right
hand side of~\thetag{ 4.7}. Since $\Omega_{ > r'} = \Omega_M \cap \{
\vert \! \vert z \vert \! \vert > r' \}$, we have $q \not \in \Omega_M$.

Reminding ${\sf R}_{ r'} \subset {\sf S}_{ r'}^{ 2n -1}$, the tangent
plane $T_p {\sf S}_{ r'}^{ 2n -1} = T_p {\sf R}_{ r'}$ divides $\C^n$
in two {\it closed}\, half-spaces, $\overline{ T}_p^+ {\sf S}_{ r'}^{
2n -1}$ exterior to $\B_{ r'}^n$ and the opposite one $\overline{
T}_p^- {\sf S}_{ r'}^{ 2n -1}$. We distinguish two (nonexclusive)
cases.

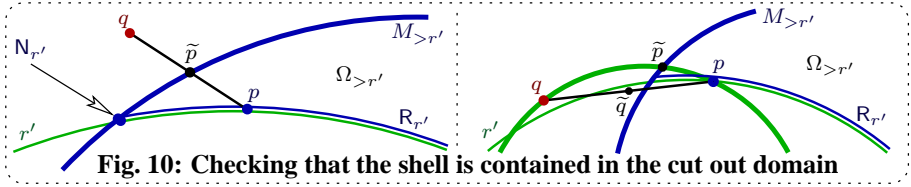
\begin{center}
\input shell-contained.pstex_t
\end{center}

Firstly, suppose that the half-line $(pq )^+$ is contained in
$\overline{ T }_p^+ {\sf S }_{ r'}^{ 2n -1}$, as in the left
figure. Since $p \in \Omega_M$ and $q \not \in \Omega_M$, there exists
at least one point $\widetilde{ p}$ of the open segment $(p, q)$ which
belongs to $M$, hence $\widetilde{ p}
\in M_{ > r }$. Then
\[
\vert\!\vert
q
-
\widetilde{p}
\vert\!\vert
<
\vert\!\vert
q-p
\vert\!\vert
<
\delta,
\] 
whence $q \in \B^n (\widetilde{ p}, \delta) \cap \{ \vert \! \vert z
\vert \! \vert > r \}$ and we deduce that $q \in \mathcal{ V
}_\delta \big( M_{ >r} \big)_{ >r}$ belongs to the right hand side
of~\thetag{ 4.7}, contradiction.

Secondly, suppose that the half-line $(p q)^+$ is contained in
$\overline{ T }_p^- {\sf S }_{ r'}^{ 2n -1}$, as in the right figure.
Let $\widetilde{ q} \in ( p, q)$ be the middle point. In the plane
passing through $0$, $p$ and $q$, consider a circle passing through
$p$ and $q$ and centered at some point close to $0$ in the open
segment $(0, \widetilde{ q })$. It has radius $< r'$ close to
$r'$. The open arc of circle between $p$ and $q$ is fully contained
in $\{ \vert \! \vert z \vert \! \vert > r' \}$.

Since $p \in \Omega_M$ and $q \not \in \Omega_M$, there exists at
least one point $\widetilde{ p}$ of the open arc of circle between $p$
and $q$ which belongs to $M$, hence $\widetilde{ p} \in M_{ > r}$.
But then $(p, q)$ is the hypothenuse of the triangle $p q \widetilde{
p}$ (remind $r' > 1$ and $\vert \! \vert q - p \vert \! \vert < \delta
< \! \! < 1$), whence $\vert \! \vert q - \widetilde{ p} \vert\! \vert
< \vert\! \vert q-p \vert\! \vert < \delta$, hence again as in the
first case, we deduce that $q \in \mathcal{ V }_\delta \big( M_{ >r}
\big)_{ >r}$, contradiction.
\endproof

If the slightly smaller radius 
\[
r ''
:=
r'
-
\eta 
= 
\widehat{r}_\kappa 
-
(k+1)
\eta
\] 
is already $< r$, we will shrink to $\eta ' := \widehat{ r}_\kappa - r
- k \eta < \eta$ the thickness of the final extensional
rind. Otherwise, in the generic case, $\widehat{ r}_\kappa - ( k + 1)
\eta$ is still $> r$. The final (exceptional) case being formally
similar, we continue the proof with $r' = \widehat{ r}_\kappa - k
\eta$ and $r'' = r' - \eta$, assuming that $r'' \geqslant r$.

Setting $r ' := \widehat{ r }_\kappa - k\eta$ in the auxiliary
Lemma~4.6, functions holomorphic in $\Omega_{ > r'} \cup \mathcal{
V}_\delta \big( M_{ > r} \big)_{ > r}$ restrict to $\text{ \rm
Shell}_{ r'}^{ r' + \delta} \big( {\sf R}_{ r'} \cup {\sf N}_{ r'}
\big)$ and then, thanks to Proposition~3.7, these restricted functions
extend holomorphically to ${\sf Rind} \big( {\sf R}_{ r'}, \eta \big)$.

\def\thelemma{4.8}\begin{lemma}
The following intersection of two open sets is connected{\rm :}
\def\theequation{4.9}\begin{equation}
{\sf Rind}
\big(
{\sf R}_{r'},\,\eta
\big)
\bigcap
\Big(
\Omega_{>r'}
\cup
\mathcal{V}_\delta
\big(
M_{>r}
\big)_{>r}
\Big).
\end{equation}
Furthermore, the union of the same two open sets contains
\def\theequation{4.10}\begin{equation}
\Omega_{>r'-\eta}
\cup
\mathcal{V}_\delta
\big(
M_{>r}
\big)_{>r}.
\end{equation}
\end{lemma}

Thus we get unique holomorphic extension to~\thetag{ 4.10} and
finally, by induction on $k$ and taking account of the final step
where $\eta$ should be shrunk appropriately, we get unique holomorphic
extension to $\Omega_{ >r } \cup \mathcal{ V }_\delta \big( M_{ >r}
\big)_{ >r}$.

The number of utilized Levi-Hartogs figures is majorated by the
product of the number of needed rinds $\sim \frac{ \widehat{ r}_\kappa
- r}{ \eta }$ times the maximal area of ${\sf R}_{ r'}$, which we
roughly majorate by the area $C \, (\widehat{ r }_\kappa )^{ 2n -1}$ of
the biggest sphere ${\sf S}_{ \widehat{ r}_\kappa }^{ 2n -1}$,
everything being divided by the area $c \, \delta^{ 2n - 1}$ covered
by a small Levi-Hartogs figure. This yields the finite number claimed
in Proposition~4.5, achieving its proof.
\endproof

\proof[Proof of Lemma~4.8] [May be skipped in a first reading] To
establish connectedness, we decompose the rind as
\[
\aligned
{\sf Rind}^+
&
:=
\big\{
(1+s)z:\
z\in{\sf R}_{r'},\
0<s<\eta/r'
\big\}
\\
{\sf Rind}^0
&
:=
{\sf R}_{r'},
\\
{\sf Rind}^-
&
:=
\big\{
(1-s)z:\
z\in{\sf R}_{r'},\
0<s<\eta/r'
\big\},
\endaligned
\]
so that ${\sf Rind} = {\sf Rind}^- \cup {\sf Rind}^0 \cup {\sf
Rind}^+$, without writing the common argument $\big( {\sf R}_{ r'},
\eta \big)$.

Obviously, the upper ${\sf Rind}^+$ is diffeomorphic to ${\sf R}_{ r'}
\times ( 0, \eta) \simeq \R^{ 2n -1} \times ( 0, \eta)$, hence is
connected. We claim that, moreover, the full ${\sf Rind}^+$ is
contained in $\Omega_{ >r'} \cup \mathcal{ V}_\delta \big( M_{ >r}
\big)_{ >r}$, whence
\def\theequation{4.11}\begin{equation}
{\sf Rind}^+
=
{\sf Rind}^+
\bigcap
\Big(
\Omega_{>r'}
\cup
\mathcal{V}_\delta
\big(
M_{>r}
\big)_{>r}
\Big).
\end{equation}

Indeed, let $q' \in {\sf Rind}^+$, hence of the form $q' = ( 1 + s)
p'$ for some $p' \in {\sf Rind}^0 = {\sf R}_{ r'}$ and some $s$ with
$0 < s < \eta / r'$. If the half-open-closed segment $(p', q']$ is
contained in $\Omega_M$, hence in $\Omega_{ >r'} = \Omega_M \cap
\big\{ \vert \! \vert z \vert \! \vert > r' \big\}$, we get for free
$q' \in \Omega_{ >r'}$.

If on the contrary, $(p', q']$ is {\it not}\, contained in $\Omega_M$,
then there exists a point $\widetilde{ q }' \in ( p', q']$ with
$\widetilde{ q}' \in M = \partial \Omega_M$, whence $\widetilde{ q}'
\in M_{ > r'} \subset M_{ >r}$ (remind $r' - \eta \geqslant r$). The
ball $\B^n ( \widetilde{ q}', \delta)$ then contains $q'$, because
$\vert \! \vert q' - \widetilde{ q }' \vert \! \vert < \vert \! \vert
q' - p' \vert \! \vert \leqslant \eta < \! \! < \delta$. This shows
$q' \in \mathcal{ V }_\delta \big( M_{ >r} \big)_{ > r}$, achieving
the claim.

\smallskip

Thus, the (upper) subpart~\thetag{ 4.11} of the intersection~\thetag{
4.9} is already connected.

\smallskip

To conclude the proof of connectedness, it suffices to show that
every point $p'$ of the remaining part
\def\theequation{4.12}\begin{equation}
\Big(
{\sf Rind}^0
\cup
{\sf Rind}^-
\Big)
\bigcap
\Big(
\Omega_{>r'}
\cup
\mathcal{V}_\delta
\big(
M_{>r}
\big)_{>r}
\Big)
\end{equation}
can be joined, by means of some appropriate continuous curve running
inside the intersection~\thetag{ 4.9}, to some point $q'$ of the
connected upper subpart~\thetag{ 4.11}. Thus, let $p'$ in~\thetag{
4.12} be arbitrary.

If $p' \in {\sf Rind}^0 \cap \big( \Omega_{ >r'} \cup \mathcal{ V
}_\delta \big( M_{>r} \big)_{ >r} \big)$, it suffices to join radially
$p'$ to $q' = (1 + s_\varepsilon ) p'$, for some $s_\varepsilon$ with
$0 < s_\varepsilon < \! \! < \eta$. Indeed, such a $q'$ then belongs to
${\sf Rind}^+ \cap \big( \Omega_{ >r'} \cup \mathcal{ V }_\delta \big(
M_{>r} \big)_{ >r} \big)$.

If $p' \in {\sf Rind}^- \cap \big( \Omega_{ >r'} \cup \mathcal{ V
}_\delta \big( M_{>r} \big)_{ >r} \big)$, then necessarily $p' \in
\mathcal{ V}_\delta \big( M_{ > r} \big)_{ > r}$, because by
definition:
\[
{\sf Rind}^-
\big(
{\sf R}_{r'},\eta
\big)
\cap 
\Omega_{>r'}
= 
\emptyset.
\] 
So there is a point $q \in M_{ >r}$ with $p' \in
\B^n (q, \delta)$.

\begin{center}
\input move-the-ball.pstex_t
\end{center}

We then distinguish two exclusive cases: either $r( q) \geqslant r'$
or $r( q) < r'$.

Firstly, assume $r (q) \geqslant r'$ (left diagram). 

If $0$, $p'$ and $q$ are aligned, we simply join $p'$ to the point $q'
:= ( 1 + s_\varepsilon ) \, \frac{ r'}{ r ( p')} \, p'$ which belongs
to ${\sf Rind }^+$. The segment $[p', q']$ is then entirely contained
in ${\sf Rind} \cap \B^n ( q, \delta)_{ > r}$, hence in~\thetag{ 4.9}.

Otherwise, in the unique plane passing through $0$, $p'$ and $q$,
consider the point $q'' := \frac{ r ( p')}{ r ( q)} \, q$, satisfying
$r ( q'') = r (p')$ and belonging to $( 0, q)$. Since $q''$ is the
orthogonal projection of $q$ onto $\overline{ \B^n ( 0, r (p'))}$, we
get $\vert \! \vert q - q'' \vert \! \vert < \vert \! \vert q - p'
\vert \! \vert < \delta$, whence $q'' \in \B^n ( q, \delta)$. The
circle of radius $r ( p')$ centered at $0$ joins $p'$ to $q''$ by
means of a small arc which is entirely contained in $\B^n ( q,
\delta)$. Denote by $\gamma : [ 0, 1] \to \B^n ( q, \delta)$ a
parametrization of this arc of circle, with $\gamma ( 0) = p'$ and
$\gamma ( 1) = q''$.

If $\gamma [ 0, 1]$ is entirely contained in ${\sf Rind}^-$, we conclude
by joining $q''$ radially to the point $q' := ( 1 + s_\varepsilon) \,
\frac{ r'}{ r( q'')}\, q''$.

If $\gamma[ 0, 1]$ is not contained in ${\sf Rind}$, let $t_1 \in ( 0,
1)$ satisfying $\gamma [0, t_1) \subset {\sf Rind}^-$ but $\gamma (
t_1) \not \in {\sf Rind}^-$. Then $\gamma ( t_1)$ belongs to $\partial
{\sf Rind}^-$ and since $r ( \gamma ( t_1)) = r ( p')$ still satisfies
$r' - \eta < r( p') < r'$, necessarily $\gamma ( t_1)$ belongs
``vertical part'' of $\partial {\sf Rind}^-$, namely to the strip
$\big\{ ( 1 - s) z : \ z \in {\sf N}_{ r'}, \ 0 \leqslant s \leqslant
\eta / r' \big\}$. Hence the point $q''' := \frac{ r'}{ r ( \gamma (
t_1))}\, \gamma ( t_1)$ belongs to ${\sf N}_{ r'}$. We now modify
$\gamma$ by constructing a curve which remains entirely inside $\B^n (
q''', \delta)_{ >r} \subset \mathcal{ V}_\delta \big( M_{ >r} \big)_{
>r}$ as follows: choose $t_2 < t_1$ very close to $t_1$, join $p'$ to
$\gamma ( t_2) \in {\sf Rind}^-$ through $\gamma$ and then $\gamma (
t_2)$ radially to the point $q' := ( 1 + s_\varepsilon) \, \frac{ r'}{
r ( \gamma (t_2))} \, \gamma (t_2) \in {\sf Rind}^+$. The resulting
curve is entirely contained in~\thetag{ 4.9}. In conclusion, we have
joined $p'$ to a suitable point $q'$, as announced.

\smallskip

Secondly, assume that $r ( q) < r'$. Consider the normalized gradient
vector field $\frac{ \nabla r_M}{ \vert \! \vert \nabla r_M \vert \!
\vert}$, defined and nowhere singular on $M \cap \big\{ \widehat{ r}_{
\kappa - 1} < \vert \! \vert z \vert \! \vert < \widehat{ r}_\kappa
\big\}$, hence on $M_{ > r} \big \backslash \{ \widehat{ p}_\kappa
\}$. For $t \in [ 0, 2\, \eta]$, denote by $t \mapsto q_t$ the integral
curve of $\frac{ \nabla r_M}{ \vert \! \vert \nabla r_M \vert \!
\vert}$ passing through $q$, satisfying $q_0 = q$, $q_t \in M$ and $r
( q_t) = r ( q) + t$. Together with its center $q$, the ball is translated
as $\B^n ( q_t, \delta)$. Accordingly, the point $p'$ is moved,
yielding a curve $p_t'$ such that $p_t'$ occupies a fixed position
with respect to the moving ball. Explicitly: $p_t' = p' + q_t' - q$.
Thanks to $r' > 1$ and $\delta < \! \! < 1$, one may check\footnote{
If the spheres ${\sf S}_{ r}^{ 2n - 1 }$ for $r$ close to $r'$ would
be hyperplanes\,\,---\,\,they almost are in comparison to $\B^n ( q_t,
\delta)$\,\,---\,\,we would have exactly $r ( p_t') = r(p') + t$,
whence $\frac{ dr ( p_t')}{ dt} = 1$.} that $\frac{ d r ( p_t')}{ dt}
\geqslant 1 - c_{ r', \delta}$, for some small positive constant $c_{
r', \delta} < \! \! 1$.

Thus, as $t$ increases, the point $p_t'$ moves away from
$0$ at speed almost equal to $1$. Since $r' - \eta < r ( p_0' ) <
r'$, we deduce that for $t = 2\, \eta$, we have $r ( p_{ 2 \eta}') > r'$,
namely $p_{ 2\eta}'$ has escaped from ${\sf Rind}^-$. Consequently,
there exists $t_1 \in ( 0, 2\eta)$ with $p_t' \in {\sf Rind}^-$ for $0
\leqslant t < t_1$ such that $p_{ t_1}' \in \partial {\sf Rind}^-$.

The boundary of ${\sf Rind}^-$ has three parts: the top ${\sf R}_{
r'}$, the bottom $\big\{ (1 - \eta /r') z : \, z \in {\sf R}_{ r'}
\big\}$ and the (closed) strip $\big\{ (1 - s) z : \, z \in {\sf N}_{
r'}, \, 0 \leqslant s \leqslant \eta / r'\big\}$. The limit point $p_{
t_1}'$ cannot belong to the bottom, since $r ( p_{ t_1}' ) 
> r ( p_0') > r' - \eta$.

Since by construction $p_t' \in \B^n ( q_t, \delta)$ with $q_t \in M_{
> r}$, we observe that $p_t' \in \mathcal{ V}_\delta \big( M_{ > r}
\big)_{ >r}$ for every $t \in [ 0, 2\, \eta]$. Consequently:
\[
p_t'
\in
{\sf Rind}^-
\bigcap
\mathcal{V}_\delta
\big(
M_{>r}
\big)_{>r},
\ \ \ \ \ 
\forall\
t\in[0,t_1).
\]

Assuming that $p_{ t_1 }' \in \partial {\sf Rind}^-$ belongs to the
top ${\sf R}_{ r'} = {\sf Rind}^0$, we may join $p_{ t_1}'$ radially
to $q' := ( 1+ s_\varepsilon) p_{ t_1 }'$. In this way, $p'$ is
joined, by means of a continuous curve running in the
intersection~\thetag{ 4.9 }, to the point $q' = (1 + s_\varepsilon)
p_{ t_1}'$ belonging to the connected upper subpart~\thetag{ 4.11}.

Finally, assume that $p_{ t_1}' \in \partial {\sf Rind}^-$ belongs to
the strip $\big\{ (1 - s) z : \, z \in {\sf N}_{ r'}, \, 0 \leqslant s
\leqslant \eta / r'\big\}$. The point $q'' := \frac{ r'}{ r ( p_{
t_1}')} \, p_{ t_1}'$ belongs to ${\sf N}_{ r'} \subset M_{ >r}$, and
we will construct a small curve running entirely inside $\B^n ( q'',
\delta)_{ >r} \subset \mathcal{ V}_\delta \big( M_{ > r} \big)_{
>r}$. Choose $t_2 \in (0, t_1)$ very close to $t_1$, join $p'$ to $p_{
t_2}' \in {\sf Rind}^-$ as above (but do not go up to $p_{ t_1}'$) and
then join $p_{ t_2}'$ radially to the point $q' := ( 1+ s_\varepsilon)
\frac{ r'}{ r ( p_{ t_2}')} \, p_{ t_2}'$, which belongs to ${\sf
Rind}^+$. The small radial segment from $p_{ t_2}'$ to $q'$ is
entirely contained in $\B^n ( q'', \delta)$ and in the full ${\sf
Rind}$. In conclusion, $p'$ is joined, by means of a continuous curve
running in the intersection~\thetag{ 4.9}, to this point $q' = ( 1+
s_\varepsilon) \frac{ r'}{ r ( p_{ t_2}')} \, p_{ t_2}'$ which belongs
to the connected upper subpart~\thetag{ 4.11}.

The proof of the connectedness of the intersection~\thetag{ 4.9}
is complete.

\medskip

We now show that the union, instead of the intersection in~\thetag{
4.9}, contains~\thetag{ 4.10}. 

Let $p' \in \Omega_{ > r' - \eta} \backslash \Omega_{ > r'}$, whence $r'
- \eta < \vert \! \vert p' \vert \! \vert \leqslant r'$. The radial
half line $\big\{ t \, p' : \, 0 < t < \infty \big\}$ emanating from
the origin and passing through $p'$ meets ${\sf S}_{ r'}^{ 2n - 1}$ at
the point $q' = \frac{ r'}{ \vert \! \vert p' \vert \! \vert }\, p'$.

If the closed segment $[ p', q']$ is contained in $\Omega_{ > r' -
\eta}$, then $q' \in \Omega_M$. Since $\vert \! \vert q' \vert \!
\vert = r'$ and since ${\sf R}_{ r'} = \Omega_M \cap \big\{ \vert \!
\vert z \vert \! \vert = r'\big\}$, we get $q' \in {\sf R}_{ r'}$,
whence $p' \in {\sf Rind} \big( {\sf R}_{ r'}, \eta \big)$.

If on the contrary, the closed segment $[ p', q']$ is {\it not}\,
contained in $\Omega_{ > r' - \eta}$, then there exists $\widetilde{
q}' \in ( p', q']$ with $\widetilde{ q}' \in M = \partial \Omega_M$,
whence $\widetilde{ q}' \in M_{ > r' - \eta} \subset M_{ > r}$. Since
$\eta < \! \! < \delta$, we deduce $p' \in \B^n ( \widetilde{ q}',
\delta)$ and finally $p' \in \mathcal{ V}_\delta \big( M_{ >r }
\big)_{ >r}$.

The proofs of Lemma~4.8 and hence also of
Proposition~4.5 are complete.
\endproof

\section*{ \S5.~Creating domains, merging 
\\
and suppressing connected components}

\subsection*{ 5.1.~Topological stability and global extensional 
geometry between regular values of $r_M$} In the preceding Section~4,
for $r$ with $\widehat{ r}_{ \kappa - 1} < r < \widehat{ r}_\kappa$,
we described the simple shape of the cut out domain $\Omega_{ >r } =
\Omega_M \cap \{ \vert \! \vert z \vert \! \vert > r \}$, just
diffeomorphic to a piece of ball. Decreasing the radius under
$\widehat{ r}_{ \kappa - 1}$, the topological picture becomes more
complex. At least for radii comprised between two singular values of
$r (z) \big\vert_M$, Morse theory assures geometrical control
together with constancy properties.

\def\thelemma{5.2}\begin{lemma}
Fix a radius $r$ satisfying $\widehat{
r}_\lambda < r < \widehat{ r}_{ \lambda + 1}$ for some $\lambda$ with
$1 \leqslant \lambda \leqslant \kappa - 1$, hence noncritical for the
distance function $r (z) \vert_M$. Then{\rm :}

\begin{itemize}

\smallskip\item[{\bf (a)}]
$T_z M + T_z {\sf S}_r^{ 2n - 1} = T_z \C^n$ at
every point $z \in M \cap {\sf S}_r^{ 2n - 1}${\rm ;}

\smallskip\item[{\bf (b)}]
the intersection $M \cap {\sf S }_r^{ 2n - 1}$ is a $\mathcal{ C
}^\infty$ compact hypersurface ${\sf N}_r \subset {\sf S}_r^{
2n -1}$ of codimension $2$ in $\C^n$, without boundary and having
finitely many connected components{\rm ;}

\smallskip\item[{\bf (c)}]
${\sf N}_{ r''}$ is diffeomorphic to ${\sf N}_{ r'}$, whenever
$\widehat{ r}_\lambda < r'' < r' < \widehat{ r}_{ \lambda
+1}${\rm ;}

\smallskip\item[{\bf (d)}]
$M_{ > r} = M \cap \{ \vert \! \vert z \vert \! \vert > r \}$ has
finitely many connected components $M_{ > r }^c$, with $1\leqslant c
\leqslant c_\lambda$, for some $c_\lambda < \infty$ which is
independent of $r${\rm ;}

\smallskip\item[{\bf (e)}]
$M_{ > r''}^c$ is diffeomorphic to $M_{ > r'}^c$, whenever $\widehat{
r}_\lambda < r'' < r' < \widehat{ r}_{ \lambda +1}$, for all $c$ with
$1 \leqslant c \leqslant c_\lambda${\rm ;}

\smallskip\item[{\bf (f)}]
$M \cap \{ r'' < \vert \! \vert z \vert \! \vert < r' \}$ is
diffeomorphic to ${\sf N}_{ r'} \times (r'', r')$, hence also to ${\sf
N}_{ r''} \times ( r'', r')$, whenever $\widehat{ r}_\lambda < r'' < r' <
\widehat{ r}_{ \lambda +1}${\rm ;}

\end{itemize}\smallskip

\end{lemma}

\proof
We summarize the known arguments of proof ({\it cf.} \cite{ mi1963}
and \cite{ hi1976}, Ch.~6).  Equivalently, {\bf (a)} says that $dr :
T_z M \to T_{ r(z)} \R$ is onto, and this holds true since by
assumption $M \cap \big\{ \widehat{ r}_\lambda < \vert \! \vert z
\vert \! \vert < \widehat{ r}_{ \lambda + 1} \big\}$ contains no
critical points of $r (z) \vert_M$.  Then {\bf (b)} follows from this
transversality {\bf (a)}.

Next, consider the Euclidean metric $( v, w ) := \sum_{ k=1}^{ 2n}\,
v_k \, w_k$ on $\C^n \simeq \R^{ 2n}$, which induces a Riemannian
metric $( \cdot , \cdot )_M$ on $M$, a nondegenerate positive bilinear
form on $TM$.  The gradient $\nabla ( r \vert_M )$ of $r(z)
\vert_M$ is the vector field on $M$ defined by requiring that $\big(
\nabla ( r\vert_M), X \big)_M = d \big( r\vert_M \big) (X)$ for all
$\mathcal{ C}^\infty$ (locally defined) vector fields $X$ on
$M$. Let ${\sf D} := 2\, {\rm Re}\, \sum_{ k= 1}^n \, z_k \,
\frac{\partial }{ \partial z_k}$ be the radial vector field which is
obviously orthogonal to spheres
and consider the orthogonal projection $X_{\sf
D}$ of ${\sf D} \vert_M$ on $TM$, a $\mathcal{ C}^\infty$ vector field
on $M$. We want to scale the gradient as
${\sf V}_{r,M} := \lambda \cdot \nabla ( r\vert_M )$ so that its
radial component is identically equal to one, namely, 
so that $\big( {\sf V}_{
r, M}, {\sf D} \big) \equiv 1$, which gives the equation:
\[
1
=
\lambda\,
\big(
\nabla(r\vert_M),\,
{\sf D}
\big)
=
\lambda\,
\big(
\nabla(r\vert_M),\,
X_{\sf D}
\big)
=
\lambda\,
\big(
\nabla(r\vert_M),\,
X_{\sf D}
\big)_M
=
\lambda\,
d\big(r\vert_M\big)(X_{\sf D}).
\]
To simply set $\lambda := \frac{ 1}{ d (r\vert_M)(X_{\sf D})}$, we
must establish that $X_{\sf D}$ cannot belong to ${\rm Ker}\, d\big(
r\vert_M \big)$ at any point $z \in M \cap \big\{ \widehat{ r}_\lambda
< \vert \! \vert z \vert \! \vert < \widehat{ r}_{ \lambda + 1}
\big\}$ of a noncritical shell.

We check this.  At such a point $z$, ${\sf D} (z)$ is not orthogonal
to $T_z M$ (otherwise $T_z M$ would coincide with $T_z {\sf S}_{ \vert
\!  \vert z \vert \! \vert}^{ 2n - 1}$), whence its orthogonal
projection $X_{\sf D} (z)$ is $\neq 0$. By definition, $({\sf D} -
X_{\sf D}) ( z)$ is orthogonal to $T_z M \ni X_{\sf D} (z)$, hence it
is orthogonal to $X_{\sf D} ( z)$ inside the $2$-dimensional plane $\Pi_z$
generated by $X_{\sf D} ( z) \neq 0$ and by ${\sf D} ( z) \neq 0$.
If, contrary to what we want, $X_{\sf D} ( z)$ would belong to ${\sf
Ker} \, d \big( r\vert_M \big) = T_z {\sf S}_{ \vert \! \vert z \vert
\! \vert}^{ 2n - 1}$, then it would be orthogonal to ${\sf D} (z)$,
and in the plane $\Pi_z$, we would have both ${\sf D} (z)$ and the
hypothenuse $\big( {\sf D} - X_{\sf D} \big) (z)$ orthogonal to
$X_{\sf D} ( z)$, which is impossible.

Thus, in spherical coordinates $( r, \vartheta_1, \dots, \vartheta_{
2n -1} )$ restricted to a noncritical shell, the $r$-component of the
$\mathcal{ C}^\infty$ scaled gradient vector field ${\sf V}_{ r, M} :=
\frac{ \nabla (r \vert_M) }{ ( \nabla (r \vert_M), \, {\sf D} )}$ is
$\equiv 1$. We deduce that the flow (wherever defined) $z_s := \exp (
s\, {\sf V}_{ r, M} \big) (z)$ simply increases the norm as $\vert \!
\vert z_s \vert \! \vert = \vert \! \vert z \vert \! \vert + s$,
whence $\exp \big( (r' - r'') {\sf V}_{ r, M} \big) (\cdot)$ induces a
diffeomorphism from ${\sf N}_{ r''}$ onto ${\sf N}_{ r'}$: this yields
{\bf (c)}.  Also, $(z'', s) \longmapsto \exp \big( (r'' + s ) {\sf
V}_{ r, M} \big) ( z '')$ gives the diffeomorphism of ${\sf N}_{ \sf
r''} \times (r' - r'')$ onto the strip $M \cap \big\{ r'' < \vert \!
\vert z \vert \! \vert < r' \big\}$, which is {\bf (f)}.

Next, the compact manifold with boundary $M_{> r} \cup {\sf N}_r$
surely has finitely many connected components, whose number is
constant for all $\widehat{ r}_\lambda < r < \widehat{ r}_{ \lambda +
1}$, because when $r$ increases or decreases, the connected components
of the slices ${\sf N}_r$ do slide smoothly in ${\sf S}_r^{ 2n - 1}$
without encountering each other: this is {\bf (d)}. Finally, {\bf (e)}
follows from {\bf (f)} and the trivial fact that the two segments
$(r'', r^0)$ and $(r', r^0)$ are diffeomorphic, whenever $\widehat{
r}_\lambda < r'' < r' < r^0 < \widehat{ r}_{ \lambda + 1}$.
\endproof

We can now state the very main technical proposition of this paper.

\def\theproposition{5.3}\begin{proposition} 
Fix a radius $r$ satisfying $\widehat{ r}_\lambda < r < \widehat{ r}_{
\lambda + 1}$ for some $\lambda$ with $1 \leqslant \lambda \leqslant
\kappa - 1$ and let $M_{ > r }^c$, $c = 1, \dots, c_\lambda$, denote
the collection of connected components of $M \cap \{ \vert \! \vert z
\vert \! \vert > r\}$. Then{\rm :}

\begin{itemize}

\smallskip

\item[{\bf (i)}]
each $M_{ > r}^c$ bounds in $\{ \vert \! \vert z \vert \! \vert > r
\}$ a {\rm unique} domain $\widetilde{ \Omega}_{ > r}^c$ which is
relatively compact in $\C^n${\rm ;}

\smallskip\item[{\bf (ii)}]
the boundary in $\C^n$ of each $\widetilde{\Omega}_{ > r}^c$,
namely{\rm :}
\[
\partial\widetilde{\Omega}_{>r}^c
=
M_{>r}^c
\cup 
{\sf N}_r^c
\cup 
\widetilde{\sf R}_r^c
\] 
consists of $M_{ >r }^c$ together with some appropriate union ${\sf
N}_r^c$ of finitely many connected components of ${\sf N}_r = M \cap
\{ \vert \! \vert z \vert \! \vert = r \}$ and with an appropriate
region $\widetilde{\sf R }_r^c \subset {\sf S }_r^{ 2n -1}$ delimited by
${\sf N }_r^c${\rm ;}

\smallskip\item[{\bf (iii)}]
two such domains $\widetilde{ \Omega }_{ > r}^{ c_1}$ and $\widetilde{
\Omega}_{ > r}^{ c_2 }$, associated to two different connected
components $M_{ >r}^{ c_1}$ and $M_{ > r}^{ c_2 }$ of $M_{ >r}$, are
either disjoint or one is contained in the other{\rm ;}

\smallskip\item[{\bf (iv)}]
for $c_1 \neq c_2$, the regions $\widetilde{ \sf R}_r^{ c_1}$ and
$\widetilde{ \sf R}_r^{ c_2}$ are either disjoint or one is contained in
the other, while their boundaries ${\sf N}_r^{ c_1}$ and 
${\sf N}_r^{ c_2}$ are
always disjoint{\rm ;}

\smallskip\item[{\bf (v)}]
for each $c = 1, \dots, c_\lambda$, every function $f$ holomorphic in
$\mathcal{ V }_\delta \big( M_{ >r} \big)_{ > r}$ has a restriction to
$\mathcal{ V }_\delta \big( M_{ >r}^c \big)_{ > r}$ which extends
holomorphically and uniquely to $\widetilde{ \Omega}_{ > r}^c$ by
means of a finite number of Levi-Hartogs figures.

\end{itemize}\smallskip

\end{proposition}

We point out that in {\bf (i)} and {\bf (ii)}, neither $\widetilde{
\Omega}_r^c$ nor $\widetilde{ \sf R}_r^c$ need be contained in our
original domain $\Omega_M$ (as it was the case in Section~4 for
$\widehat{ r}_{ \kappa - 1} < r < \widehat{ r}_\kappa$): this is why
we introduced a widetilde notation. We refer to the middle Figure~1
for an illustration. Similarly, neither $\widetilde{ \Omega }_r^c$
nor $\widetilde{ \sf R}_r^c$ need be contained in $\C^n \big
\backslash \overline{ \Omega}_M$: they both may intersect $\Omega_M$
{\it and}\, $\C^n \big \backslash \overline{ \Omega}_M$. Also, the
number of connected components of ${\sf N }_r^c$ is $\geqslant$ that
of $\widetilde{ \sf R }_r^c$ and may be $>$, as illustrated below.

\begin{center}
\begin{picture}(0,0)%
\includegraphics{various-topological-situations.pstex}%
\end{picture}%
\setlength{\unitlength}{4144sp}%
\begingroup\makeatletter\ifx\SetFigFont\undefined
\def\x#1#2#3#4#5#6#7\relax{\def\x{#1#2#3#4#5#6}}%
\expandafter\x\fmtname xxxxxx\relax \def\y{splain}%
\ifx\x\y   
\gdef\SetFigFont#1#2#3{%
  \ifnum #1<17\tiny\else \ifnum #1<20\small\else
  \ifnum #1<24\normalsize\else \ifnum #1<29\large\else
  \ifnum #1<34\Large\else \ifnum #1<41\LARGE\else
     \huge\fi\fi\fi\fi\fi\fi
  \csname #3\endcsname}%
\else
\gdef\SetFigFont#1#2#3{\begingroup
  \count@#1\relax \ifnum 25<\count@\count@25\fi
  \def\x{\endgroup\@setsize\SetFigFont{#2pt}}%
  \expandafter\x
    \csname \romannumeral\the\count@ pt\expandafter\endcsname
    \csname @\romannumeral\the\count@ pt\endcsname
  \csname #3\endcsname}%
\fi
\fi\endgroup
\begin{picture}(5424,1374)(439,-1873)
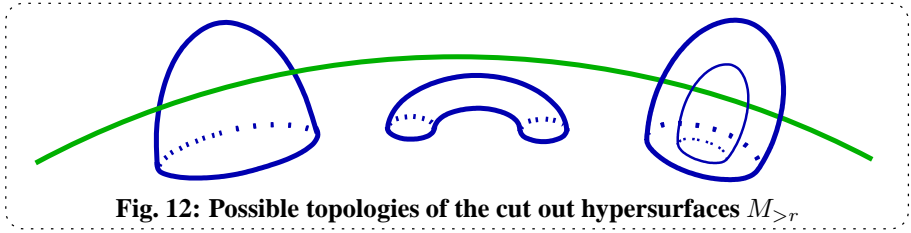
\put(1109,-1789){\makebox(0,0)[lb]{\smash{\SetFigFont{10}{12.0}{rm}{\color[rgb]{0,0,0}{\bf Fig.~12: Possible topologies of the cut out  hypersurfaces $M_{ > r}$}}%
}}}
\end{picture}

\end{center}

As a direct application, we may achieve the proof of our principal
result.

\def\thetheorem{5.4}\begin{theorem}
Under the precise assumptions of Theorem~2.7, holomorphic functions in
$\mathcal{ V}_\delta ( M)$ do extend holomorphically and
uniquely to $\Omega_M$ by means
of a finite number of Levi-Hartogs figures{\rm :}
\[
\forall\ 
f\in\mathcal{O}\big(\mathcal{V}_\delta(M)\big)
\ \ \ \ \
\exists\
F\in\mathcal{O}
\big(
\Omega_M\cup
\mathcal{V}_\delta(M)
\big)
\ \ \ 
\text{\rm s.t.}
\ \ 
F\big\vert_{\mathcal{V}_\delta(M)}
=f.
\]
\end{theorem}

\proof
In the main Proposition~5.3, we choose $r = \widehat{ r}_1 +
\varepsilon$ (where $\varepsilon > 0$ satisfies $\varepsilon < \! \! <
\delta$) very close to the last, smallest singular radius. Then $M_{ >
r}$ has a single connected component, $M_{ > r}$ itself, and it simply
bounds $\big( \Omega_M \big)_{ >r}$. The remainder part of $M$, namely
$M \cap \big\{ \vert \! \vert z \vert \! \vert \leqslant \widehat{
r}_1 + \varepsilon \big\}$ is diffeomorphic to a very small closed $(2
n - 1)$-dimensional spherical cap and is entirely contained in
$\mathcal{ V}_\delta (M)$. 

Fix an arbitrary function $f \in \mathcal{ O} \big(\mathcal{ V}_\delta
( M) \big)$ and restrict it to $\mathcal{ V}_\delta \big( M_{ >r}
\big)_{ >r}$. Thanks to the proposition, $f$ extend holomorphically
and uniquely to $\big( \Omega_M \big)_{ > r}$ by means of a finite
number of Levi-Hartogs figures. Since
\[
\mathcal{V}_\delta(M)
\bigcap
\Big(
\mathcal{V}_\delta 
\big(M_{>r}\big)_{>r}
\cup
\big(
\Omega_M\big)_{>r}
\Big)
\]
is easily seen to be 
connected, we get a globally defined extended function which is
holomorphic in
\[
\mathcal{V}_\delta(M)
\bigcup
\Big(
\mathcal{V}_\delta 
\big(M_{>r}\big)_{>r}
\cup
\big(\Omega_M\big)_{>r}
\Big)
=
\mathcal{V}_\delta(M)
\cup
\Omega_M.
\]
This completes the proof.
\endproof

\proof[Proof of Proposition~5.3]
In {\bf (i)}, let us
check the uniqueness of a relatively compact $\widetilde{
\Omega}_{ > r}^c$. Since $M_{ >r}^c$ inherits an
orientation from $M$, the complement $\big\{ \vert \! \vert z \vert \!
\vert > r \big\} \big \backslash M_{ > r}^c$ has at most 2 connected
components. As $M \Subset \C^n$ is bounded, at least one component
contains the points at infinity, hence there can remain at most one
component of $\big\{ \vert \! \vert z \vert \! \vert > r\big\} \big
\backslash M_{ > r}^c$ that is relatively compact in $\C^n$.

\smallskip

If $r$ satisfies $\widehat{ r}_{ \kappa - 1} < r < \widehat{
r}_\kappa$, Proposition~4.5 already completes the proof.

\smallskip

Assume therefore that $r$ satisfies $\widehat{ r}_\mu < r <
\widehat{ r}_{ \mu + 1}$, for some $\mu \in \N$ with $1 \leqslant \mu
\leqslant \kappa - 1$. For every $\lambda$ with $2 \leqslant
\lambda \leqslant \kappa - 1$, it will be convenient to flank each
singular radius $\widehat{ r}_\lambda$ by the following two very close
nonsingular radii
\def\theequation{5.5}\begin{equation}
\boxed{\widehat{r}_\lambda^-
:=
\widehat{r}_\lambda
-
\eta/2}
\quad\quad
\text{\rm and}
\quad\quad
\boxed{
\widehat{r}_\lambda^+
:=
\widehat{r}_\lambda
+
\eta/2},
\end{equation}
with $\eta$ being the same uniform thickness of extensional rinds as
before. We fix once for all an arbitrary function $f$ holomorphic in
$\mathcal{ V}_\delta \big( M_{ >r } \big)_{ >r}$. Letting $\lambda$ be
arbitrary with $\mu \leqslant \lambda \leqslant \kappa -1$, the logic
of the proof shows up two topologically distinct phenomena that we
overview.

\medskip\noindent
{\sf A: Filling domains through regular radii intervals.} Assume
that at the regular radius $\widehat{ r}_{ \lambda + 1}^- = \widehat{
r}_{ \lambda + 1} - \frac{ \eta}{ 2}$, all domains $\widetilde{
\Omega}_{ > \widehat{ r}_{ \lambda + 1 }^-}^c$, $c = 1, \dots,
c_\lambda$, as well as the corresponding holomorphic extensions, have
been constructed. Then prolong the domains (without topological
change) as $\widetilde{ \Omega}_{ > \widehat{ r}_\lambda^+}^c$, $c = 1,
\dots, c_\lambda$, up to $\widehat{ r}_\lambda^+ = \widehat{
r}_\lambda + \frac{ \eta}{ 2}$ and fill in the conquered
territory by means of a finite number of Levi-Hartogs figures.

\medskip\noindent
{\sf B: Jumping across singular radii and changing the domains.}
Restarting at $\widehat{ r}_\lambda^+$ with the domains $\widetilde{
\Omega}_{ > \widehat{ r}_\lambda^+}^c$, $c = 1, \dots, c_\lambda$,
distinguish three cases as follows. Remind from \S2.3 that $M$ is
represented by $v = \sum_{ 1 \leqslant j \leqslant k_\lambda} \, x_j^2
- \sum_{ 1 \leqslant j \leqslant 2n - k_\lambda - 1} y_j^2$ in
suitable coordinates $(x, y, v)$ centered at $\widehat{ p}_\lambda$,
where $k_\lambda$ is the {\sl Morse coindex} of $r ( z)
\vert_M$ at $\widehat{ p}_\lambda$.

\begin{itemize}
\smallskip\item[{\bf (I)}]
Firstly, assume $k_\lambda = 0$, namely $z \mapsto r ( z) \vert_M$ has
a local maximum at $\widehat{ p}_\lambda$, or inversely, 
assume $k_\lambda =
2n - 1$, namely $z \mapsto r ( z) \vert_M$ has a local minimum at
$\widehat{ p}_\lambda$. This is the easiest case, the only one in
which new domains can be born or die, locally.

\smallskip\item[{\bf (II)}]
Secondly, assume $k_\lambda = 1$. This is the most delicate case,
because in a small neighborhood of $\widehat{ p}_\lambda$, the cut out
hypersurface $M_{ > \widehat{ r}_\lambda^+}$ has exactly 2 connected
components, so that two different enclosed domains $\widetilde{
\Omega}_{ > \widehat{ r}_\lambda^+}^{ c_1}$ and $\widetilde{ \Omega}_{ >
\widehat{ r}_\lambda^+}^{ c_2}$ can meet here; it may also occur that
the two parts near $\widehat{ p}_\lambda$ belong to the {\it same}\,
domain, i.e. that $c_2 = c_1$. While descending down to $\widehat{
r}_\lambda^-$, we must analyze the way how the two (maybe the single)
component(s) merge. Three subcases will be distinguished, one of
which showing a crucial trick of {\it subtracting}\, one growing
component from a larger one which also grows
(right Figure~1).

\smallskip\item[{\bf (III)}]
Thirdly, assume that $2 \leqslant k_\lambda \leqslant 2n - 2$. In all
these cases, locally in a neighborhood of $\widehat{ p}_\lambda$, the
cut out hypersurface $M_{ > \widehat{ r}_\lambda^+}$ has exactly 1
connected component and the way how the corresponding single enclosed
domain $\widetilde{ \Omega}_{ > \widehat{ r}_\lambda^+}^c$ grows will
be topologically constant.

\end{itemize}\smallskip

\smallskip\noindent

Reasoning by induction on $\lambda$ and applying the filling processes
{\sf A} and {\sf B}, we then descend progressively inside deeper
spherical shells, checking all properties of Proposition~5.3. When
approaching the bottom radius $r$ of Proposition~5.3, it will suffice
to shortcut {\sf A} or {\sf B} appropriately in order to complete the
proof.

\subsection*{ 5.6.~Filling domains through regular radii 
intervals}
Recall that $\widehat{ r}_\mu < r < \widehat{ r}_{ \mu + 1}$, let
$\lambda$ with $\mu \leqslant \lambda \leqslant \kappa - 1$ and
consider the regular radius interval $\big[ \widehat{ r}_\lambda^+ ,
\widehat{ r}_{ \lambda + 1}^- \big]$. We suppose first that $r
\leqslant \widehat{ r}_\lambda^+$, so that we may descend inside the
whole spherical shell $\big\{ \widehat{ r}_\lambda^+ < \vert \! \vert
z \vert \! \vert \leqslant \widehat{ r}_{ \lambda + 1}^- \big\}$.
Afterwards, we explain how we stop in the case where $\lambda = \mu$
and $\widehat{ r}_\mu^+ < r < \widehat{ r}_{ \mu + 1}^-$.

\smallskip

By descending induction on $\lambda$ through
{\sf A} and {\sf B}, we may assume that at the
superlevel set
$(\cdot)_{ > \widehat{ r}_{ \lambda + 1 }^-}$, the domains $\widetilde{
\Omega}_{ > \widehat{ r}_{ \lambda + 1}^-}^c$ enclosed by $M_{ >
\widehat{ r}_{ \lambda + 1}^-}^c$ for $1 \leqslant c \leqslant
c_\lambda$ have been constructed and that each restriction $f_{
\widehat{ r}_{ \lambda + 1}^-}^c$ of $f \in \mathcal{ O} \big(
\mathcal{ V}_\delta \big( M_{ > r} \big)_{ > r} \big)$ to $\mathcal{
V}_\delta \big( M_{ > \widehat{ r}_{ \lambda + 1}^-}^c \big)_{ >
\widehat{ r}_{ \lambda + 1}^-}$ extends holomorphically and uniquely
to the domain
\def\theequation{5.7}\begin{equation}
\widetilde{\Omega}_{ 
>\widehat{ r}_{\lambda +1}^-}^c
\bigcup 
\mathcal{V}_\delta 
\big(M_{>
\widehat{r}_{\lambda+1}^-}^c\big)_{> 
\widehat{r}_{\lambda+1}^-}.
\end{equation}

For every radius $r'$ with $\widehat{ r}_\lambda^+ \leqslant r' <
\widehat{ r}_{ \lambda + 1}^-$, the cut out hypersurface $M_{ > r'} =
\bigcup_{ 1 \leqslant c \leqslant c_\lambda} \, M_{ > r'}^c$ has the
same number of connected components, each $M_{ > r'}^c$ is
diffeomorphic to $M_{ > \widehat{ r}_{ \lambda + 1}^-}^c$ and the
difference $M_{ > r'}^c \big \backslash M_{ > \widehat{ r}_{ \lambda
+1}^-}^c$ is diffeomorphic to $N_{ \widehat{ r}_{ \lambda +1}^-}^c
\times \big( r', \widehat{ r}_{ \lambda +1}^- \big]$. Furthermore,
each prolongation $\widetilde{ \Omega}_{ > r'}^c$ of $\widetilde{
\Omega}_{ > \widehat{ r}_{ \lambda + 1}^-}^c$ is obviously defined
just by adding the tube domain surrounded by $M_{ > r'}^c \big
\backslash M_{ \widehat{ r}_{ \lambda + 1}^- }^c$. Then each $N_{ r'
}^c = \partial \widetilde{ \sf R}_{ r' }^c$ has finitely many
connected components $N_{ r'}^{ c, j}$, with $1 \leqslant j \leqslant
j_{ \lambda, c}$, where $j_{ \lambda, c}$ is independent of $r'$.

\begin{center}
\input legs.pstex_t
\end{center}

Since $f$ was defined in $\mathcal{ V}_\delta \big( M_{ > r} \big)_{ >
r}$ and since $r \leqslant \widehat{ r}_\lambda^+$, we claim that each
restriction $f_{ \widehat{ r}_{ \lambda + 1}^-}^c$ may be extended
holomorphically and uniquely to
\def\theequation{5.8}\begin{equation}
\widetilde{\Omega}_{>\widehat{r}_{\lambda+1}^-}^c
\bigcup
\mathcal{V}_\delta
\big(
M_{>\widehat{r}_\lambda^+}^c
\big)_{>\widehat{r}_\lambda^+}.
\end{equation}
Indeed, to the original domain of definition~\thetag{ 5.7} of $f_{
\widehat{ r}_{ \lambda + 1}^-}^c$ which was contained in $\big\{ \vert
\! \vert z \vert \! \vert > \widehat{ r}_{ \lambda + 1}^- \big\}$, we
add in the enlarged domain~\thetag{ 5.8} a finite number $j_{ \lambda,
c}$ of tubular domains around the connected components of $M_{ > r'}^c
\big \backslash M_{ > \widehat{ r}_{ \lambda +1}^-}^c$. Because
$\delta$ was chosen so small that $\mathcal{ V}_\delta ( M)$ is a
small tubular neighborhood of $M$, 
and because $f \in \mathcal{ O} \big(
\mathcal{ V}_\delta
\big( M_{ > r} \big)_{ > r} \big)$
is uniquely defined, we get a unique
extension, still denoted by $f_{ \widehat{ r}_{ \lambda +
1}^-}^c$, to~\thetag{ 5.8}.

We can now apply the same reasoning as in Proposition~4.5, which
consists of progressive holomorphic extension by means of thin
rinds. Reproducing the proof of Lemma~4.6 (with changes of notation
only), we get for every radius $r'$ with $\widehat{ r}_\lambda^+ < r'
\leqslant \widehat{ r}_{ \lambda + 1}^-$ that
\def\theequation{5.9}\begin{equation}
\text{\rm Shell}_{r'}^{r'+\delta}
\big(
\widetilde{\sf R}_{r'}^c
\cup 
{\sf N}_{r'}^c
\big)
\ \
\text{\rm is contained in}
\ \
\widetilde{\Omega}_{>r'}^c
\bigcup
\mathcal{V}_{\delta}
\big(
M_{>\widehat{r}_\lambda^+}^c
\big)_{>\widehat{r}_\lambda^+}.
\end{equation}
Similarly, reproducing the proof of Lemma~4.8 yields the connectedness
of
\[
{\sf Rind}
\big(
{\sf R}_{r'}^c,\,\eta
\big)
\bigcap
\Big(
\widetilde{\Omega}_{>r'}^c
\cup
\mathcal{V}_{\delta}
\big(
M_{>\widehat{r}_\lambda^+}^c
\big)_{>\widehat{r}_\lambda^+}
\Big),
\]
and furthermore, this yields that the union, 
instead of the intersection, contains
\[
\widetilde{\Omega}_{>r'-\eta}^c
\bigcup
\mathcal{V}_{\delta}
\big(
M_{>\widehat{r}_\lambda^+}^c
\big)_{>\widehat{r}_\lambda^+},
\]
whenever $r'- \eta$ is still $\geqslant \widehat{ r}_\lambda^+$
(otherwise, shrink conveniently the thickness of the last extensional
rind, as in the proof of Proposition~4.5). Thus, by piling up $\frac{
\widehat{ r}_{ \lambda + 1}^- - \widehat{ r}_\lambda^+}{ \eta}$ rinds
and by using a finite number $\leqslant C \, \big( \frac{ \widehat{
r}_\kappa}{ \delta} \big)^{ 2n -1} \, \Big[ \frac{ \widehat{ r}_{
\lambda + 1}^- - \widehat{ r}_\lambda^+}{ \eta} \Big]$ of Levi-Hartogs
figures, we get unique holomorphic extension to
\def\theequation{5.10}\begin{equation}
\mathcal{V}_\delta
\big(
M_{>\widehat{r}_\lambda^+}^c
\big)_{>\widehat{r}_\lambda^+}
\bigcup
\widetilde{\Omega}_{>\widehat{r}_\lambda^+}^c.
\end{equation}

Finally, if $r$ satisfies $\widehat{ r}_\mu^+ < r < \widehat{ r}_{ \mu
+ 1}^-$, descending from $( \cdot)_{ > \widehat{ r}_{ \mu + 1}^-}$
with $\lambda = \mu$ as above, we just stop the construction of rinds
to $( \cdot)_{ >r}$ by shrinking appropriately the thickness of the
last extensional rind.

The property {\bf (iii)} that enclosed domains $\widetilde{ \Omega}_{
> r }^c$ are either disjoint or one is contained in the other remains
stable as $r$ decreases through the whole nonsingular interval $\big(
\widehat{ r}_\lambda, \widehat{ r}_{ \lambda + 1} \big)$, because
their (moving) boundaries always remain disjoint, so that property
{\bf (iv)} is also simultaneously transmitted to lower regular radii.
This completes {\sf A}.

\subsection*{ 5.11.~Localizing (pseudo)cubes at Morse points}
We now study {\sf B}. Recall that $\widehat{ r}_\mu < r < \widehat{
r}_{ \mu + 1}$, let $\lambda$ with $\mu \leqslant \lambda \leqslant
\kappa - 1$ and suppose that $r \leqslant \widehat{ r}_{ \lambda }^-$,
so that starting from $(\cdot)_{ > \widehat{ r}_\lambda^+}$, we may
(and we must) continue the Hartogs-Levi filling inside the whole thin
spherical shell $\big\{ \widehat{ r}_{ \lambda }^- < \vert \! \vert z
\vert \! \vert \leqslant \widehat{ r}_\lambda^+ \big\}$. Similarly as
above, the way how we should stop the process in the case where
$\lambda = \mu$ and $\widehat{ r}_\mu < r < \widehat{ r}_{ \mu }^+$
is obvious.

\smallskip

By descending induction on $\lambda$ through {\sf A} and {\sf B}, we
may assume that at $\widehat{ r}_\lambda^+$, the domains $\widetilde{
\Omega}_{ > \widehat{ r}_\lambda^+ }^c$ enclosed by $M_{ > \widehat{
r}_\lambda^+}^c$ for $1 \leqslant c \leqslant c_\lambda$ have been
constructed and that each restriction $f_{ \widehat{ r}_\lambda^+ }^c$
of $f \in \mathcal{ O} \big( \mathcal{ V}_\delta \big( M_{ > r}^c
\big)_{ > r} \big)$ to $\mathcal{ V}_\delta \big( M_{ > \widehat{
r}_\lambda^+}^c \big)_{ > \widehat{ r}_\lambda^+}$ extends
holomorphically to the domain~\thetag{ 5.10} of the previous paragraph.

By an elementary analysis of the Morse normalizing quadric, we will
see that in some small (pseudo)cube centered at $\widehat{
p}_\lambda$, there passes in most cases only one component $M_{ >
\widehat{ r}_\lambda^+ }^c$, while in a single exceptional case, there
can pass two (at most) different
connected components $M_{ > \widehat{ r}_\lambda^+ }^{
c_1}$ and $M_{ > \widehat{ r}_\lambda^+ }^{ c_2}$. We will consider
only this single (or these two) component(s), because the other
components do pass regularly and without topological change accross
$\widehat{ p}_\lambda$, hence are filled in 
by Levi-Hartogs figures exactly as in {\sf A}.

\smallskip

Shrinking the $\delta_1$ of Theorem~2.7 if necessary (remind $0 <
\delta \leqslant \delta_1$), we may assume that the Morse normalizing
coordinates $\big( v, x_1, \dots, x_{ k_\lambda}, y_1, \dots, y_{ 2n -
1 - k_\lambda} \big)$ near $\widehat{ p}_\lambda$ are defined in the
ball $\B^n ( \widehat{ p}_\lambda, \delta_1 )$ and that the map
\[
z
\longmapsto
\big(
v(z),x(z),y(z)
\big), 
\ \ \ \ \ \ \ \ \ \
\B^n ( \widehat{ p}_\lambda, \delta_1 )
\longrightarrow
\R^{2n}
\]
is close in $\mathcal{ C }^1$ norm to its differential at $\widehat{
p}_\lambda$, so that it is almost not distorting. Then $\delta_1$ shall
not be shrunk anymore.

Because in the estimates of the (finite) number of Levi-Hartogs
figures, $\eta$ only appears as a denominator in a factor $\frac{ r' -
r''}{ \eta}$ ({\it cf.} Proposition~4.5), it is allowed to work with
extensional rinds of smaller universal positive thickness, at the cost
of spending a number of pushed analytic discs that is greater, of
course, but still finite. If necessary, we shrink $\eta >0$ to insure
that $\eta^{ 1 / 2} < \! \! < \delta$. Then $\eta$ will not be shrunk
anymore.

Thanks to these preliminaries, we may define a convenient
(pseudo)cube centered at $\widehat{ p }_\lambda$ by
\def\theequation{5.12}\begin{equation}
{\sf C}_\eta
:=
\Big\{
z\in
\B^n(\widehat{p}_\lambda,
\delta_1):\
\vert 
v(z)
\vert
<
\eta,\ \
\vert\!\vert
x(z)
\vert\!\vert
<
2\,\eta^{1/2},\ \ 
\vert\!\vert
y(z)
\vert\!\vert
<
2\,\eta^{1/2}
\Big\}.
\end{equation}
It then follows that ${\sf C}_\eta$ is properly contained in
$\mathcal{ V}_\delta ( M)$ and is relatively small. Reminding that
$v(z) = r (z) - r ( \widehat{ p}_\lambda)$, the radial thickness of
${\sf C}_\eta$ is equal to $2 \eta$, twice the difference $\widehat{
r}_\lambda^+ - \widehat{ r}_\lambda^- = \eta$. We draw a diagram
assuming $k_\lambda = 2n -1$ (see only the left one).

\begin{center}
\input C-eta.pstex_t
\end{center}

\subsection*{ 5.13.~Topology of horizontal super-level sets in the
complement of quadrics} Simultaneously to the proof, we provide an
auxiliary elementary study. Let $n\in \N$ with $n\geqslant 2$, let $k
\in \N$ with $0 \leqslant k \leqslant 2n - 1$, let $x = (x_1, \dots,
x_k) \in \R^k$, let $y = (y_1, \dots, y_{ 2n - 1 - k} ) \in \R^{ 2n -
1 - k}$, let $v \in \R$, and in $\R^{ 2n }$ equipped with the
coordinates $(x, y, v)$, consider the quadric of equation
\def\theequation{5.14}\begin{equation}
v
=
\sum_{1\leqslant j\leqslant k}\,
x_j^2
-
\sum_{1\leqslant j\leqslant 2n-1-k}\,
y_j^2,
\end{equation}
which we will denote by ${\sf Q}_k$. The coordinate $v$ playing the
r\^ole of $r ( z) - r (\widehat{ p}_\lambda)$ near a singular radius
$\widehat{ r}_\lambda$ having Morse coindex $k_\lambda$, we want to
understand how the topology of the super-level sets
\[
\big\{v>\varepsilon\big\} 
\cap
\big(\R^{2n}\backslash 
{\sf Q}_k
\big)
\]
(which relate to the possible domains $\widetilde{ \Omega}_{ > r}^c$
for $r$ close to $\widehat{ r}_\lambda$) do change when the parameter
$\varepsilon$ descends from a small positive value to a small negative
value.

\begin{center}
\begin{picture}(0,0)%
\includegraphics{filling-cap.pstex}%
\end{picture}%
\setlength{\unitlength}{4144sp}%
\begingroup\makeatletter\ifx\SetFigFont\undefined
\def\x#1#2#3#4#5#6#7\relax{\def\x{#1#2#3#4#5#6}}%
\expandafter\x\fmtname xxxxxx\relax \def\y{splain}%
\ifx\x\y   
\gdef\SetFigFont#1#2#3{%
  \ifnum #1<17\tiny\else \ifnum #1<20\small\else
  \ifnum #1<24\normalsize\else \ifnum #1<29\large\else
  \ifnum #1<34\Large\else \ifnum #1<41\LARGE\else
     \huge\fi\fi\fi\fi\fi\fi
  \csname #3\endcsname}%
\else
\gdef\SetFigFont#1#2#3{\begingroup
  \count@#1\relax \ifnum 25<\count@\count@25\fi
  \def\x{\endgroup\@setsize\SetFigFont{#2pt}}%
  \expandafter\x
    \csname \romannumeral\the\count@ pt\expandafter\endcsname
    \csname @\romannumeral\the\count@ pt\endcsname
  \csname #3\endcsname}%
\fi
\fi\endgroup
\begin{picture}(5424,1284)(439,-1063)
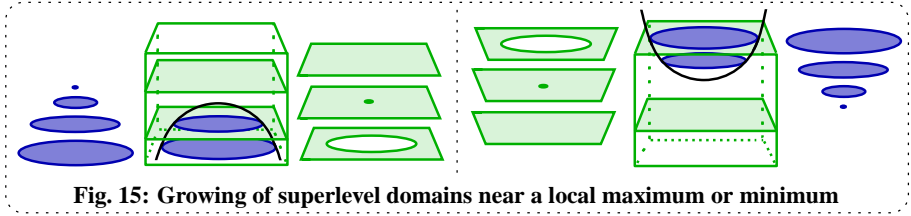
\put(856,-990){\makebox(0,0)[lb]{\smash{\SetFigFont{9}{10.8}{rm}{\color[rgb]{0,0,0}{\bf Fig.~15: Growing of superlevel domains near a local maximum or minimum}}%
}}}
\end{picture}

\end{center}

In the case $k= 0$ (left figure) the quadric looks like a spherical
cap, its complement $\R^{ 2 n } \big \backslash {\sf Q}_0$ having
exactly two connected components. For positive values of
$\varepsilon$, there is only one (green) super-level component $\{ v >
\varepsilon \} \cap \big( \R^{ 2 n} \big \backslash { \sf Q}_0
\big)$. As $\varepsilon$ becomes negative, this component grows
regularly, allowing a newly created hole to widen inside the slices
$\{ v = \varepsilon \}$. The (blue) holes then pile up to constitute a
newly created, local component $M_{ > \widehat{ r}_\lambda^- }^c$.

The (reverse) case $k = 2n - 1$ exhibits the
local end of some component $M_{ > \widehat{ r }_\lambda^- }^c$. In
a while, we will see that there is a salient topological difference
between the two remaining (less obvious) cases 
$2 \leqslant k \leqslant 2n -
2$ and $k=1$, the exceptional one. Before pursuing, we conclude the
proof of {\sf B} in case $\widehat{ p}_\lambda$ is a local maximum or
minimum.

\smallskip

We assume $k_\lambda = 2n - 1$, the case $k_\lambda = 0$ being already
considered (essentially completely) in Section~4. Observe that $M_{ >
\widehat{ r }_\lambda^+ } \cap {\sf C}_\eta$ is diffeomorphic to
${\sf S}^{ 2n - 2} \times (c/2, c)$, hence connected. Thus, let $M_{ >
\widehat{ r }_\lambda^+ }^c$ denote the single component entering
${\sf C}_\eta$. By descending induction through {\sf A} and {\sf B},
$M_{ > \widehat{ r }_\lambda^+ }^c$ bounds a relatively compact
domain of holomorphic extension $\widetilde{ \Omega}_{ > \widehat{
r}_\lambda^+}^c$, with $\partial \widetilde{ \Omega}_{ > \widehat{
r}_\lambda^+}^c = M_{ > \widehat{ r}_\lambda^+}^c \cup {\sf N}_{
\widehat{ r}_\lambda^+}^c \cup \widetilde{ \sf R}_{ \widehat{
r}_\lambda^+}^c$, as in property {\bf (ii)} of Proposition~5.3, all
the other properties also holding true on $( \cdot )_{ > \widehat{
r}_\lambda^+}$ . Denote by $\widetilde{ \sf R}_{ \widehat{
r}_\lambda^+}^{ c, k}$, $1 \leqslant k \leqslant k_{ \lambda, c}$, the
connected components of $\widetilde{ \sf R}_{ \widehat{
r}_\lambda^+}^c$ and by ${ \sf N}_{ \widehat{ r}_\lambda^+}^{ c, j}$,
$1 \leqslant j \leqslant j_{ \lambda, c}$, with $j_{ \lambda, c}
\geqslant k_{ \lambda, c}$, the components of ${\sf N}_{ \widehat{
r}_\lambda^+}^c$.

\begin{center}
\input fill-k-2n-1.pstex_t
\end{center}

\noindent
We do the numbering so that ${\sf C}_\eta$ encloses the first (small)
${\sf N}_{ \widehat{ r}_\lambda^+}^{c, 1}$, which is diffeomorphic to
a small $(2n - 2)$-dimensional sphere. Also, we number so that the
boundary of $\widetilde{ \sf R}_{ \widehat{ r}_\lambda^+ }^{c, 1}$ in
${\sf S}_{ \widehat{ r}_\lambda^+ }^{ 2n - 1}$ contains ${\sf N}_{
\widehat{ r }_\lambda^+ }^{c, 1}$, whence $\widetilde{ \sf R}_{
\widehat{ r }_\lambda^+ }^{c, 1}$ meets ${\sf C}_\eta$. We do not
draw ${\sf C}_\eta$.

\smallskip

Observe that, by means of extensional rinds that are symmetric around
the other components $\widetilde{ \sf R}_{ \widehat{ r }_\lambda^+ }^{
c, 2}, \dots, \widetilde{ \sf R}_{ \widehat{ r }_\lambda^+ }^{ c, k_{
\lambda, c} }$, we may achieve the Hartogs-Levi filling exactly as in
{\sf A}, because $r(z) \vert_M$ is regular in $\mathcal{ V}_\delta
\big( {\sf N}_{ \widehat{ r }_\lambda^+ }^{ c, j} \big)$, for every
$j$ such that ${\sf N}_{ \widehat{ r }_\lambda^+ }^{ c, j}$ is
contained in the boundary of each of these other components. Hence it
remains only to discuss what is happening in a neighborhood of the
single component $\widetilde{ \sf R}_{ \widehat{ r}_\lambda^+}^{ c,
1}$, and especially near $\widehat{ p }_\lambda$.

For the disposition of $\widetilde{ \Omega}_{ > \widehat{
r}_\lambda^+}^c \cap {\sf C}_\eta$, or equivalently of $\widetilde{
\sf R}_{ \widehat{ r}_\lambda^+}^{ c, 1} \cap {\sf C}_\eta$, two cases
occur. Let $\big( v, x_1, \dots, x_{ 2n - 1} \big)$ be the Morse
coordinates centered at $\widehat{ p}_\lambda$.

\begin{itemize}

\smallskip\item[{\sf (a)}]
As illustrated by the left figure above, $\widetilde{ \Omega}_{ >
\widehat{ r }_\lambda^+ }^c \cap {\sf C }_\eta$ consists of the
space\footnote{ Sets written ``$\{ \cdot \}$'' here are understood to
be subsets of ${\sf C}_\eta$.} lying above $\big\{ v = \eta / 2\big\}$
and above $\big\{ v = x_1^2 + \cdots + x_{ 2n - 1}^2 \big\}$, a
cap-shaped space which is clearly connected; the region $\widetilde{
\sf R}_{ \widehat{ r}_\lambda^+ }^{c, 1}$ is then diffeomorphic to a
small $(2n - 1)$-dimensional ball.

\smallskip\item[{\sf (b)}]
As illustrated by the right figure above, $\widetilde{ \Omega}_{ >
\widehat{ r }_\lambda^+}^c \cap {\sf C}_\eta$ consists of the space
lying above $\big\{ v = \eta / 2 \big\}$ but below $\big\{ v = x_1^2 +
\cdots + x_{ 2n - 1}^2 \big\}$~; the dimension of ${\sf S}_{ \widehat{
r}_\lambda^+}^{ 2n -1}$ being $\geqslant 3$, the region $\widetilde{
R}_{ \widehat{ r }_\lambda^+ }^{ c, 1} \cap {\sf C}_\eta$ is
connected, a fact that a one-dimensional diagram cannot show
adequately; then $\widetilde{ \Omega }_{ > \widehat{ r }_\lambda^+}^c
\cap {\sf C}_\eta$ is also connected.

\end{itemize}\smallskip

In case {\sf (a)}, near $\widehat{ p }_\lambda$, a piece of
$\widetilde{ \Omega}_{ > \widehat{ r }_\lambda^+ }^c$ ends up while
descending to the lower super-level set $( \cdot )_{ > \widehat{ r
}_\lambda^- }$. We do not use any extensional rind there, we just
observe that unique holomorphic extension is got for free in
\[
\Big[
\mathcal{V}_\delta
\big(
M_{>\widehat{r}_\lambda^-}^c
\big)_{>\widehat{r}_\lambda^-}
\Big]
\cap
{\sf C}_\eta,
\]
since this domain is fully contained in $\mathcal{ V }_\delta \big(
M_{ > r} \big)_{ > r}$.

\smallskip

In case {\sf (b)}, we apply Hartogs Levi
extension to ${\sf Rind} \big(
\widetilde{ \sf R }_{ \widehat{ r }_\lambda^+}^{ c, 1}, \eta \big)$
and we get unique holomorphic extension from~\thetag{ 5.10} to
\[
\Big[
\mathcal{V}_\delta
\big(
M_{>\widehat{r}_\lambda^-}^c
\big)_{>\widehat{r}_\lambda^-}
\Big]
\bigcup
{\sf Rind}
\big(
\widetilde{\sf R}_{ 
\widehat{r}_\lambda^+}^{c,1},\eta
\big).
\]
The union of this open set together with~\thetag{ 5.10} contains a
unique well defined domain $\widetilde{ \Omega}_{ > \widehat{ r
}_\lambda^- }^c$ with the property that the passage from $\widetilde{
\sf R}_{ > \widehat{ r }_\lambda^+ }^{ c, 1}$ to $\widetilde{ \sf R}_{ >
\widehat{ r}_\lambda^-}^{ c, 1}$ fills a hole, as illustrated by the right
diagram above, whence ${\sf N}_{ > \widehat{ r }_\lambda^-}^c$ has one
less connected component, because the $( 2n -2 )$-sphere ${\sf N}_{ >
\widehat{ r}_\lambda + \varepsilon }^{ c, 1}$ drops when $\varepsilon
< 0$.

The properties that two different domains $\widetilde{ \Omega}_{ >
\widehat{ r}_\lambda^+ }^{ c_1}$ and $\widetilde{ \Omega}_{ >
\widehat{ r}_\lambda^+ }^{ c_2}$ are either disjoint or one is
contained in the other is easily seen to be inherited by $\widetilde{
\Omega}_{ > \widehat{ r}_\lambda^- }^{ c_1}$ and $\widetilde{
\Omega}_{ > \widehat{ r}_\lambda^- }^{ c_2}$: it suffices to
distinguish two cases: $c_2 \neq c$ and $c_1 \neq c$, or $c_2 \neq c$
and $c_1 = c$; to look at {\sf (a)} or {\sf (b)} and then to conclude.

\smallskip

The proof of {\sf B} in case $k_\lambda = 2n - 1$ is complete. The
case $k_\lambda = 0$ is similar: two subcases {\sf 
(a')}\,\,---\,\,reverse {\sf (a)}\,\,---\,\,and {\sf
(b')}\,\,---\,\,reverse {\sf (b)}\,\,---\,\, then appear; subcase {\sf
(a')} exhibits the birth of a new component (blue left Figure~15), as
already fully studied in Section~4 while subcase {\sf (b')} (green
left Figure~15) shows that an external component descends regularly
as do clouds around a hill.

\subsection*{5.15.~The regular cases $2 \leqslant k_\lambda \leqslant
2n - 2$}
Let $k$ with $2 \leqslant k \leqslant 2n - 2$ and consider
the quadric ${\sf Q}_k$ of~\thetag{ 5.14}. We
claim that
${\sf Q}_k \cap \big\{ v > \varepsilon \big\}$
has exactly one connected component for
every $\varepsilon >0$.
Indeed, ${\sf Q}_k \cap \big\{ v > \varepsilon \big\}$ can
be represented as
\[
\bigcup_{y_1,\dots,y_{2n-k-1}}\
\bigcup_{\varepsilon'>\varepsilon}\
\big\{
x_1^2
+\cdots+
x_k^2
=
\varepsilon'
+
y_1^2
+\cdots+
y_{2n-1-k}^2
\big\}.
\]
Since $\varepsilon '$ is always positive, we hence have a smoothly
parameterized family of $(k-1)$-dimensional spheres that are all
connected. Consequently, the union is also connected, as claimed.

To view the topology more adequately, in the case $n=2$, we draw a
short movie consisting of the 3-dimensional slices $\big\{ v =
\varepsilon ' \big\} \cap \big( \R^{2n} \big \backslash {\sf Q}_k
\big)$, where $\varepsilon' = \frac{ 2}{ 3}\, \eta, \ \frac{ 1}{ 2}\,
\eta, \ 0, \ - \frac{ 1}{ 2}\, \eta$. To conceptualize (in case
$n=2$) the super-level sets
\[
\big\{
v 
> 
\varepsilon
\big\} 
\cap 
\big(
\R^{2n}
\big\backslash 
{\sf Q}_k
\big)
=
\bigcup_{\varepsilon'>\varepsilon}\,
\big\{
v 
=
\varepsilon'
\big\} 
\cap 
\big(
\R^{2n}
\big\backslash 
{\sf Q}_k
\big),
\]
it suffices to pile up
intuitively the images of the corresponding movie.

\begin{center}
\begin{picture}(0,0)%
\includegraphics{2-k-2n-2.pstex}%
\end{picture}%
\setlength{\unitlength}{4144sp}%
\begingroup\makeatletter\ifx\SetFigFont\undefined
\def\x#1#2#3#4#5#6#7\relax{\def\x{#1#2#3#4#5#6}}%
\expandafter\x\fmtname xxxxxx\relax \def\y{splain}%
\ifx\x\y   
\gdef\SetFigFont#1#2#3{%
  \ifnum #1<17\tiny\else \ifnum #1<20\small\else
  \ifnum #1<24\normalsize\else \ifnum #1<29\large\else
  \ifnum #1<34\Large\else \ifnum #1<41\LARGE\else
     \huge\fi\fi\fi\fi\fi\fi
  \csname #3\endcsname}%
\else
\gdef\SetFigFont#1#2#3{\begingroup
  \count@#1\relax \ifnum 25<\count@\count@25\fi
  \def\x{\endgroup\@setsize\SetFigFont{#2pt}}%
  \expandafter\x
    \csname \romannumeral\the\count@ pt\expandafter\endcsname
    \csname @\romannumeral\the\count@ pt\endcsname
  \csname #3\endcsname}%
\fi
\fi\endgroup
\begin{picture}(5424,1286)(439,-1108)
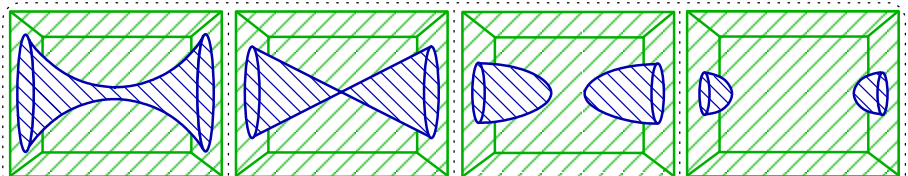
\put(514,-1031){\makebox(0,0)[lb]{\smash{\SetFigFont{9}{10.8}{rm}{\color[rgb]{0,0,0}{\bf Fig.~17: Sliced view of the growing of the two possible domains in case $2\leqslant k_\lambda \leqslant 2n - 2$}}%
}}}
\end{picture}

\end{center}

So let $M_{ > \widehat{ r }_\lambda^+ }^c$ be the single connected
component of $M \cap \big\{ \vert \! \vert z \vert \! \vert >
\widehat{ r}_\lambda^+ \big\}$ that enters ${\sf C}_\eta$. The
corresponding domain $\widetilde{ \Omega}_{ > \widehat{ r}_\lambda^+
}^c$ can be located from one or the other side. Its prolongation up to
the deeper sublevel set $(\cdot)_{ > \widehat{ r}_\lambda^-}$ (viewed
only inside ${\sf C}_\eta$) consists of piling up the (blue) small
symmetric regions or the (green) surrounding regions drawn above.

We do the numbering so that ${\sf N}_{ \widehat{ r }_\lambda^+}^{ c,
1}$ enters ${\sf C}_\eta$, being a (connected) hyperboloid as drawn in
the first picture of Figure~17 and so that $\widetilde{ \sf R}_{
\widehat{ r }_\lambda^+ }^{ c, 1}$ enters ${\sf C}_\eta$ as one
(connected, blue or green) side of this hyperboloid. As previously in
the two cases $k_\lambda = 0$ and $k_\lambda = 2n - 1$, the
Hartogs-Levi filling goes through exactly as in the regular case {\bf
A} for all other $\widetilde{ \sf R}_{ \widehat{ r}_\lambda^+ }^{ c,
2}, \dots, \widetilde{ \sf R}_{ \widehat{ r}_\lambda^+ }^{ c,
k_{\lambda, c}}$. Next, by putting finitely many Levi-Hartogs
figures in ${\sf Rind} \big( \widetilde{ \sf
R}_{\widehat{ r }_\lambda^+}^{ c, 1}, \eta \big)$ 
we get holomorphic extension from the
domain~\thetag{ 5.10} to
\[
\Big[
\mathcal{V}_\delta
\big(
M_{>\widehat{r}_\lambda^-}^c
\big)_{>\widehat{r}_\lambda^-}
\Big]
\bigcup
{\sf Rind}
\big(
\widetilde{\sf R}_{ 
\widehat{r}_\lambda^+}^{c,1},\eta
\big).
\]
The intersection of~\thetag{ 5.10}
with this open set is connected because 
$\widetilde{ \sf R}_{\widehat{ r}_\lambda^+}^{ c, 1}$
is connected, and the union of both 
contains a 
well defined
domain $\widetilde{ \Omega}_{ > \widehat{ r}_\lambda^-}^c$
obtained by adding the (blue or green) slices of Figure~17.


\section*{\S6.~The exceptional case $k_\lambda = 1$} 

\subsection*{6.1.~Illustration}
To begin with the most delicate case, we draw a 3-dimensional diagram
showing a saddle-like $M$ localized in a (pseudo)cube ${\sf C }_\eta$
centered at $\widehat{ p }_\lambda$.

\begin{center}
\input saddle.pstex_t
\end{center}

For every $\varepsilon$ satisfying $0 < \varepsilon < \eta$, there are
two connected components $M_{ > \varepsilon}^-$ and $M_{ >
\varepsilon}^+$ of $M_{ > \widehat{ r }_\lambda + \varepsilon} \cap
{\sf C }_\eta$, namely the two upper tips of the saddle, defined in
equations by
\[
M_{>\varepsilon}^\pm
:=
\big\{
v
=
x^2
-
y_1^2
-\cdots-
y_{2n-2}^2
\big\}
\cap
\big\{
\pm
x
>
0
\big\}
\cap
\big\{
v
>
\varepsilon
\big\}.
\]
With $\varepsilon = \frac{ 1}{ 2 }\, \eta$, we are simply looking at
$M_{ > \widehat{ r }_\lambda^+ } \cap {\sf C }_\eta$. By descending
induction through {\sf A} and {\sf B}, we are given two domains of
holomorphic extension $\widetilde{ \Omega}_{ > \widehat{ r }_\lambda^+
}^{ c_1}$ and $\widetilde{ \Omega}_{ > \widehat{ r }_\lambda^+ }^{
c_2}$ whose boundary contains $M_{ > \eta / 2}^-$ and $M_{ > \eta /
2}^+$, respectively.

Firstly, we assume that $c_2 \neq c_1$. Since each one of the two
pieces of hypersurfaces $M_{ > \eta / 2}^-$ and $M_{ > \eta / 2}^+$
has two sides, there are $2 \times 2 = 4$ subcases to be considered
for the relative disposition of $\Omega_{ > \eta / 2}^- := \widetilde{
\Omega}_{ > \widehat{ r }_\lambda^+ }^{ c_1} \cap {\sf C}_\eta$ and of
$\Omega_{ > \eta /2}^+ := \widetilde{ \Omega}_{ > \widehat{ r
}_\lambda^+ }^{ c_2} \cap {\sf C}_\eta$, with $c_2 \neq c_1$.

\begin{itemize}

\smallskip\item[{\sf (a)}]
$\Omega_{ > \eta / 2}^-$ \big(resp. $\Omega_{ > \eta / 2}^+$\big)
consists of the space lying above the hyperplane $\{ v = \eta / 2 \}$
and below the left (resp. right) tip of the saddle, namely in
equations:
\[
\Omega_{>\eta/2}^\pm
=
\big\{
v
>
\eta/2
\big\}
\bigcap
\big\{
\pm
x
>
0
\big\}
\bigcap
\big\{
v
<
x^2
-
y_1^2
-
\cdots
-
y_{2n-2}^2
\big\}.
\]

\smallskip\item[{\sf (b)}]
$\Omega_{ > \eta / 2}^-$ is the small nose as in {\sf (a)} but
$\Omega_{ > \eta / 2}^+$ consists of the other side, {\it i.e.} of the
(rather bigger) space lying inside $\big\{ v > \eta / 2 \big\}$ left
to $M_{ > \eta / 2}^+$, namely in equations:
\[
\Omega_{>\eta/2}^+
=
\big\{
v>\eta/2
\big\}
\Big\backslash
\Big(
\big\{
x>0
\big\}
\bigcap
\big\{
v
\leqslant
x^2
-
y_1^2
-\cdots-
y_{2n-2}^2
\big\}
\Big).
\]

\smallskip\item[{\sf (c)}]
Symetrically to {\sf (b)}, $\Omega_{ > \eta / 2}^+$ is the small nose as
in {\sf (a)} but 
\[
\Omega_{>\eta/2}^-
=
\big\{
v>\eta/2
\big\}
\Big\backslash
\Big(
\big\{
x<0
\big\}
\bigcap
\big\{
v
\leqslant
x^2
-
y_1^2
-\cdots-
y_{2n-2}^2
\big\}
\Big).
\]

\smallskip\item[{\sf (d)}]
Finally, $\Omega_{ > \eta / 2}^-$ is as in {\sf (c)} and $\Omega_{ >
\eta / 2}^+$ is as in {\sf (b)}.

\end{itemize}\smallskip

The last subcase {\sf (d)} cannot occur, because it is ruled out by
property {\bf (iii)} of Proposition~5.3, which holds on the
super-level set $( \cdot)_{ > \widehat{ r}_\lambda^+}$ by the inductive
assumption.

\smallskip

Secondly, we assume that $c_2 = c_1$. Then there can occur a subcase
{\sf (a')} very similar to {\sf (a)}, in which $c_2 = c_1$, so that
$\Omega_{ > \eta / 2}^-$ and $\Omega_{ > \eta / 2}^+$ belong to the
{\it same}\, enclosed relatively compact domain. But with $c_2 = c_1$, no
subcase similar to {\sf (b)}\,\,---\,\,or to {\sf (c)}\,\,---\,\,can occur,
because $M_{ > \eta / 2}^- \subset \partial \Omega_{ > \eta /
2}^+$\,\,---\,\,or $M_{ \eta / 2}^+ \subset \partial \Omega_{ > \eta /
2}^-$\,\,---\,\,would then bound the {\it same} relatively compact
domain from its both sides, but we already know from the beginning of
the proof, that one side at least must always contain the points at
infinity.

Finally, with $c_2 = c_1 = c$, there remains the following last subcase
(unseen previously).

\begin{itemize}

\smallskip\item[{\sf (e)}]
$\Omega_{ > \eta / 2} := \widetilde{ \Omega}_{ > \widehat{
r}_\lambda^+ }^c \cap {\sf C}_\eta$ consists of the space lying above
$\big\{ v = \eta / 2 \big\}$ and above the saddle, namely
\[
\Omega_{>\eta/2}
=
\big\{
v
>
\eta/2
\big\}
\bigcap
\big\{
v
>
x^2
-
y_1^2
-\cdots-
y_{2n-2}^2
\big\}.
\]

\end{itemize}\smallskip

\smallskip

As $M = \partial \Omega_M$ lies in $\C^n$ with $n \geqslant 2$, whence
$2n - 2 \geqslant 2$, there is at least one dimension of $y \in \R^{
2n - 2}$ which is missing in the left figure above. To view the
topology more adequately, coming back to the abstract quadric ${\sf
Q}_1$ and assuming $n = 2$, we plan to draw a short movie consisting
of the 3-dimensional slices $\big\{ v = \varepsilon ' \big\} \cap
\big( \R^{2n} \big \backslash {\sf Q}_1 \big)$, where $\varepsilon ' =
\frac{ 2}{ 3}\, \eta, \ \frac{ 1}{ 2}\, \eta, \ 0, \ - \frac{ 1}{ 2}\,
\eta$.

Recall that we are interested in the connected components of the
super-level sets
\[
\big\{
v 
> 
\varepsilon
\big\} 
\cap 
\big(
\R^{2n}
\big\backslash 
{\sf Q}_1
\big)
=
\bigcup_{\varepsilon'>\varepsilon}\,
\big\{
v 
=
\varepsilon'
\big\} 
\cap 
\big(
\R^{2n}
\big\backslash 
{\sf Q}_1
\big).
\]
As suggested by this sliced union, to conceptualize these
4-dimensional (in case $n=2$) super-level sets, it suffices to pile up
intuitively the images of the corresponding movie.

\begin{center}
\begin{picture}(0,0)%
\includegraphics{i-k-1.pstex}%
\end{picture}%
\setlength{\unitlength}{4144sp}%
\begingroup\makeatletter\ifx\SetFigFont\undefined
\def\x#1#2#3#4#5#6#7\relax{\def\x{#1#2#3#4#5#6}}%
\expandafter\x\fmtname xxxxxx\relax \def\y{splain}%
\ifx\x\y   
\gdef\SetFigFont#1#2#3{%
  \ifnum #1<17\tiny\else \ifnum #1<20\small\else
  \ifnum #1<24\normalsize\else \ifnum #1<29\large\else
  \ifnum #1<34\Large\else \ifnum #1<41\LARGE\else
     \huge\fi\fi\fi\fi\fi\fi
  \csname #3\endcsname}%
\else
\gdef\SetFigFont#1#2#3{\begingroup
  \count@#1\relax \ifnum 25<\count@\count@25\fi
  \def\x{\endgroup\@setsize\SetFigFont{#2pt}}%
  \expandafter\x
    \csname \romannumeral\the\count@ pt\expandafter\endcsname
    \csname @\romannumeral\the\count@ pt\endcsname
  \csname #3\endcsname}%
\fi
\fi\endgroup
\begin{picture}(5424,1286)(439,-1108)
\put(2174,-73){\makebox(0,0)[lb]{\smash{\SetFigFont{9}{10.8}{rm}{\color[rgb]{0,0,0}$\big\{ v = \frac{1}{2}\,\eta \big\}$}%
}}}
\put(3571,-124){\makebox(0,0)[lb]{\smash{\SetFigFont{9}{10.8}{rm}{\color[rgb]{0,0,0}$\big\{ v =0\big\}$}%
}}}
\put(4952,-776){\makebox(0,0)[lb]{\smash{\SetFigFont{9}{10.8}{rm}{\color[rgb]{0,0,0}{\bf merge}}%
}}}
\put(4804, -4){\makebox(0,0)[lb]{\smash{\SetFigFont{9}{10.8}{rm}{\color[rgb]{0,0,0}$\big\{ v =-\frac{1}{2} \eta \big\}$}%
}}}
\put(810,-269){\makebox(0,0)[lb]{\smash{\SetFigFont{9}{10.8}{rm}{\color[rgb]{0,0,0}$\big\{ v = \frac{2}{3}\,\eta \big\}$}%
}}}
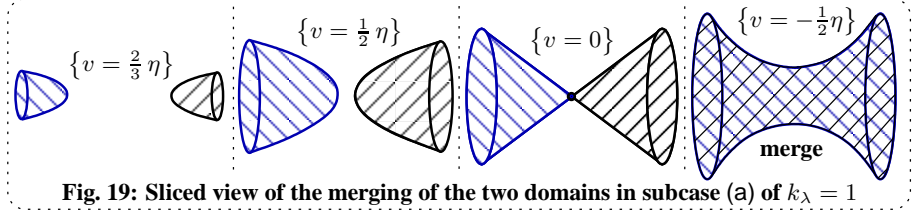
\put(775,-1042){\makebox(0,0)[lb]{\smash{\SetFigFont{9}{10.8}{rm}{\color[rgb]{0,0,0}{\bf Fig.~19: Sliced view of the merging of the two domains in subcase {\sf (a)} of $k_\lambda=1$}}%
}}}
\end{picture}

\end{center}

Here, the second picture shows $\widetilde{ \sf R}_{ \widehat{
r}_\lambda^+}^{ c_1} \cap {\sf C}_\eta$ (in blue, to the left)
together with $\widetilde{ \sf R}_{ \widehat{ r}_\lambda^+}^{ c_1}
\cap {\sf C}_\eta$ (in black, to the right). Then the third picture
shows how the two components do touch and the fourth one shows how they
should be merged as $\varepsilon ' = - \frac{ 1}{ 2}\, \eta$ becomes
negative. The complete discussion follows in a while.

We next offer the movie of {\sf (b)}, the movie of {\sf (c)} being
obtained from it just by a reflection across the hyperplane $\{ x = 0
\}$.

\begin{center}
\begin{picture}(0,0)%
\includegraphics{ii-k-1.pstex}%
\end{picture}%
\setlength{\unitlength}{4144sp}%
\begingroup\makeatletter\ifx\SetFigFont\undefined
\def\x#1#2#3#4#5#6#7\relax{\def\x{#1#2#3#4#5#6}}%
\expandafter\x\fmtname xxxxxx\relax \def\y{splain}%
\ifx\x\y   
\gdef\SetFigFont#1#2#3{%
  \ifnum #1<17\tiny\else \ifnum #1<20\small\else
  \ifnum #1<24\normalsize\else \ifnum #1<29\large\else
  \ifnum #1<34\Large\else \ifnum #1<41\LARGE\else
     \huge\fi\fi\fi\fi\fi\fi
  \csname #3\endcsname}%
\else
\gdef\SetFigFont#1#2#3{\begingroup
  \count@#1\relax \ifnum 25<\count@\count@25\fi
  \def\x{\endgroup\@setsize\SetFigFont{#2pt}}%
  \expandafter\x
    \csname \romannumeral\the\count@ pt\expandafter\endcsname
    \csname @\romannumeral\the\count@ pt\endcsname
  \csname #3\endcsname}%
\fi
\fi\endgroup
\begin{picture}(5424,1286)(439,-1108)
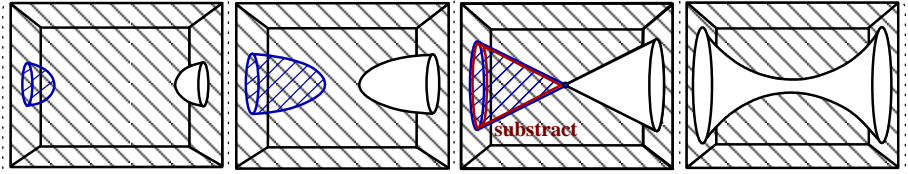
\put(664,-1042){\makebox(0,0)[lb]{\smash{\SetFigFont{9}{10.8}{rm}{\color[rgb]{0,0,0}{\bf Fig.~20: Sliced view of the substraction of the left domain in subcase {\sf (b)} of $k_\lambda=1$}}%
}}}
\put(3392,-689){\makebox(0,0)[lb]{\smash{\SetFigFont{8}{9.6}{rm}{\color[rgb]{.69,0,0}\red{\bf substract}}%
}}}
\end{picture}

\end{center}

Here again, the second picture shows $\widetilde{ \sf R}_{ \widehat{
r}_\lambda^+}^{ c_1} \cap {\sf C}_\eta$ (in blue, to the left)
together with $\widetilde{ \sf R}_{ \widehat{ r}_\lambda^+}^{ c_2}
\cap {\sf C}_\eta$ (the large (black) region, containing the small (blue)
one). Then the third picture, namely the slice $\varepsilon ' = 0$,
shows a not allowed situation: the left cone does bound {\it two}\,
regions from its {\it two}\, sides, contrary to the {\it a priori}\,
unique relatively compact domain $\widetilde{ \Omega}_{ \widehat{
r}_\lambda}^{ c_1} \subset \big\{ \vert \! \vert z \vert \! \vert >
\widehat{ r}_\lambda \big\}$ we are seeking to construct, when
starting from $\widetilde{ \Omega}_{ \widehat{ r}_\lambda^+}^{ c_1}$.
The trick is then to {\it suppress}\, the (blue)
small slice, or equivalently
to subtract it from the (black) large slice which contains it. Then the
black winning slice continues to grow up to $\big\{ v = - \eta / 2
\big\}$ (fourth picture). The complete discussion follows in a while.

Finally, here is the (simpler) movie of {\bf (e)}.

\begin{center}
\begin{picture}(0,0)%
\includegraphics{iii-k-1.pstex}%
\end{picture}%
\setlength{\unitlength}{4144sp}%
\begingroup\makeatletter\ifx\SetFigFont\undefined
\def\x#1#2#3#4#5#6#7\relax{\def\x{#1#2#3#4#5#6}}%
\expandafter\x\fmtname xxxxxx\relax \def\y{splain}%
\ifx\x\y   
\gdef\SetFigFont#1#2#3{%
  \ifnum #1<17\tiny\else \ifnum #1<20\small\else
  \ifnum #1<24\normalsize\else \ifnum #1<29\large\else
  \ifnum #1<34\Large\else \ifnum #1<41\LARGE\else
     \huge\fi\fi\fi\fi\fi\fi
  \csname #3\endcsname}%
\else
\gdef\SetFigFont#1#2#3{\begingroup
  \count@#1\relax \ifnum 25<\count@\count@25\fi
  \def\x{\endgroup\@setsize\SetFigFont{#2pt}}%
  \expandafter\x
    \csname \romannumeral\the\count@ pt\expandafter\endcsname
    \csname @\romannumeral\the\count@ pt\endcsname
  \csname #3\endcsname}%
\fi
\fi\endgroup
\begin{picture}(5424,1286)(439,-1108)
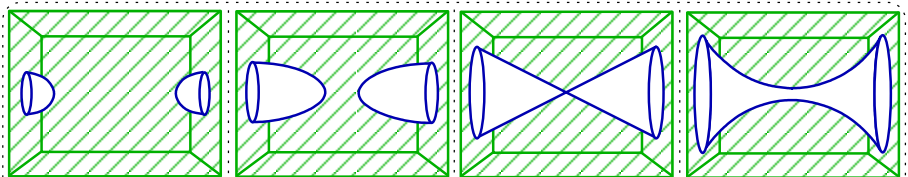
\put(693,-1042){\makebox(0,0)[lb]{\smash{\SetFigFont{9}{10.8}{rm}{\color[rgb]{0,0,0}{\bf Fig.~21: Sliced view of the growing of the external domain in subcase {\sf (e)} of $k_\lambda=1$}}%
}}}
\end{picture}

\end{center}

\subsection*{ 6.2.~Jumping across the singular radius:
merging process} Assuming $k_\lambda = 1$, we can now complete {\sf
B} in subcase {\sf (a)}, postponing subcase {\sf (a')}. We
look at Figures~17 and~18.

Let $M_{ > \widehat{ r}_\lambda^+}^{ c_1} \cap {\sf C}_\eta$ and $M_{
> \widehat{ r}_\lambda^+}^{ c_2} \cap {\sf C }_\eta$ be the two
``nose'' components of $M_{ > \widehat{ r}_\lambda^+}$ entering ${\sf
C}_\eta$. Here, $c_2 \neq c_1$. By descending
induction through {\sf A} and {\sf B}, $M_{ > \widehat{ r }_\lambda^+
}^{ c_1}$ and $M_{ > \widehat{ r }_\lambda^+ }^{ c_2}$ bound some two
relatively compact domains of holomorphic extension $\widetilde{
\Omega}_{ > \widehat{ r}_\lambda^+}^{ c_1}$ and $\widetilde{ \Omega}_{
> \widehat{ r}_\lambda^+}^{ c_2}$ with $\partial \widetilde{ \Omega}_{
> \widehat{ r}_\lambda^+}^{ c_1} = M_{ > \widehat{ r}_\lambda^+ }^{
c_1} \cup {\sf N}_{ \widehat{ r}_\lambda^+}^{ c_1} \cup \widetilde{
\sf R}_{ \widehat{ r}_\lambda^+ }^{ c_1}$ and $\partial \widetilde{
\Omega}_{ > \widehat{ r}_\lambda^+ }^{ c_2} = M_{ > \widehat{
r}_\lambda^+}^{ c_2} \cup {\sf N}_{ \widehat{ r}_\lambda^+}^{ c_2}
\cup \widetilde{ \sf R}_{ \widehat{ r}_\lambda^+ }^{ c_2}$ as in
property {\bf (ii)} of Proposition~5.3, all the other properties also
holding true on $( \cdot )_{ > \widehat{ r }_\lambda^+}$.

We remind that the other domains $\widetilde{ \Omega}_{ > \widehat{
r}_\lambda^+}^c$ for $c \neq c_1$ and $c \neq c_2$
with $1 \leqslant c \leqslant c_\lambda$ do pass regularly
through $\widehat{ r}_\lambda$ up to $( \cdot)_{ > \widehat{
r}_\lambda^-}$, thanks to {\sf A}.

For $i=1, 2$, denote by $\widetilde{ \sf R}_{ \widehat{
r}_\lambda^+}^{ c_i, k}$, $1 \leqslant k \leqslant k_{ \lambda, c_i}$,
the connected components of $\widetilde{ \sf R}_{ \widehat{
r}_\lambda^+}^{c_i}$ and by ${ \sf N}_{ \widehat{ r}_\lambda^+}^{ c_i,
j}$, $1 \leqslant j \leqslant j_{ \lambda, c_i}$, with $j_{ \lambda,
c_i} \geqslant k_{ \lambda, c_i}$, the components of ${\sf N}_{
\widehat{ r}_\lambda^+}^{ c_i}$. We do the numbering so that
$\widetilde{ \sf R}_{ \widehat{ r}_\lambda^+}^{ c_1, 1}$
(resp. $\widetilde{ \sf R}_{ \widehat{ r}_\lambda^+}^{ c_2, 1}$)
enters ${\sf C}_\eta$ to the left (resp. right), together with ${ \sf
N}_{ \widehat{ r}_\lambda^+}^{ c_1, 1}$ (resp. ${ \sf N}_{ \widehat{
r}_\lambda^+}^{ c_2, 1}$), as illustrated by
Figure~17.

As in the case $k_\lambda = 2n -1$, for $i=1, 2$, by means of
extensional rinds that are symmetric around the other components
$\widetilde{ \sf R}_{ \widehat{ r }_\lambda^+ }^{ c_i, 2}, \dots,
\widetilde{ \sf R}_{ \widehat{ r }_\lambda^+ }^{ c_i, k_{ \lambda,
c_i} }$, we may achieve the Hartogs-Levi filling exactly as in {\sf
A}, because $r(z) \vert_M$ is regular in $\mathcal{ V}_\delta \big(
{\sf N}_{ \widehat{ r }_\lambda^+ }^{ c_i, j} \big)$, for every $j$
such that ${\sf N}_{ \widehat{ r }_\lambda^+ }^{ c_i, j}$ is contained
in the boundary of each of these other components. Hence it remains
only to discuss what is happening in a neighborhood of the two
components $\widetilde{ \sf R}_{ \widehat{ r}_\lambda^+ }^{ c_i, 1}$,
$i=1, 2$, and especially near the saddle point $\widehat{ p
}_\lambda$.

\smallskip

While descending from $\widehat{ r}_\lambda^+$ to $\widehat{
r}_\lambda^-$, the two regions $\widetilde{ \sf R}_{ \widehat{
r}_\lambda^+}^{ c_1, 1} \subset {\sf S}_{ \widehat{ r}_\lambda^+}^{ 2n
-1}$ and $\widetilde{ \sf R}_{ \widehat{ r}_\lambda^+}^{ c_2, 1}
\subset {\sf S}_{ \widehat{ r}_\lambda^+}^{ 2n -1}$ do merge as a
single connected region contained in ${\sf S}_{ \widehat{
r}_\lambda^-}^{ 2n -1}$ that we will denote by $\widetilde{
\sf R}_{ \widehat{
r}_\lambda^-}^*$, {\it see}\, the right Figure~17. In Morse theory
(\cite{ mi1963, hi1976}), one speaks of {\sl attaching a one-cell},
since in the merging process, the two regions are essentially joined
by means of a (thickened) segment directed along the $x$-axis. It
follows that the two hypersurfaces $M_{ > \widehat{ r}_\lambda^+}^{
c_1}$ and $M_{ > \widehat{ r}_\lambda^+}^{ c_2}$ do merge as a
connected hypersurface $M_{ > \widehat{ r}_\lambda^-}^*$ containing
them, and furthermore, that the two domains $\widetilde{ \Omega}_{ >
\widehat{ r}_\lambda^+}^{ c_1}$ and $\widetilde{ \Omega}_{ > \widehat{
r}_\lambda^+}^{ c_2}$ do prolong uniquely up to the slightly deeper
super-level set $( \cdot)_{ > \widehat{ r }_\lambda^- }$, merging as a
uniquely defined domain $\widetilde{ \Omega}_{ > \widehat{ r
}_\lambda^- }^*$ which is relatively compact
in $\C^n$ and which contains
$\widetilde{\sf R }_{\widehat{ r
}_\lambda^+}^*$ in its boundary $\partial \widetilde{ \Omega}_{
>\widehat{r }_\lambda^+}^*$.

As $c_2 \neq c_1$, the new number of
domains in the interval $(
\widehat{ r}_{ \lambda - 1}, \widehat{ r}_\lambda)$ is lowered by a
unit, {\it i.e.} $c_{ \lambda - 1} = c_\lambda - 1$
(if $c_2 = c_1$ as in {\sf (a')}, 
the number would not change, {\it i.e.} $c_{ \lambda - 1} = c_\lambda$).

For $i=1, 2$, let $f_{ \widehat{ r}_\lambda^+}^{ c_i}$ denote the
restriction of $f \in \mathcal{ O} \big( \mathcal{ V}_\delta \big( M_{
> r} \big)_{ > r} \big)$ to $\mathcal{ V}_{ \delta} \big( M_{ >
\widehat{ r}_\lambda^+}^{ c_i} \big)_{ > \widehat{ r}_\lambda^+}$. By
descending induction through {\sf A} and {\sf B}, $f_{ \widehat{
r}_\lambda^+}^{ c_i}$ extends holomorphically and uniquely to
$\widetilde{ \Omega}_{ > \widehat{ r}_\lambda^+}^{ c_i}$. 
Then both functions do extend holomorphically and uniquely to
\[
\mathcal{V}_\delta 
\big(M_{>\widehat{r}_\lambda^-}^*
\big)_{>\widehat{ r}_\lambda^-}, 
\]
since they coincide with $f$ near $\widehat{ p}_\lambda$. 
We then introduce the two extensional rinds ${\sf Rind} \big(
\widetilde{ \sf R}_{ \widehat{ r}_\lambda^+ }^{ c_i}, \eta \big)$,
drawn in the right Figure~17. Two applications of Proposition~3.7
together with a geometrically clear connectedness property yield
unique holomorphic extension
to 
\[
{\sf Rind}
\big(
\widetilde{
\sf R}_{ 
\widehat{r}_\lambda^+ }^{c_1},\eta 
\big)
\bigcup
{\sf Rind}
\big(
\widetilde{
\sf R}_{
\widehat{r}_\lambda^+ }^{c_2},\eta 
\big)
\bigcup
\mathcal{V}_\delta 
\big(M_{>\widehat{r}_\lambda^-}^*
\big)_{\widehat{ r}_\lambda^-}
\bigcup
\widetilde{\Omega}_{>\widehat{r}_\lambda^+}^{c_1}
\bigcup
\widetilde{\Omega}_{>\widehat{r}_\lambda^+}^{c_2}.
\]
In sum, we have got unique holomorphic extension to
\[
\mathcal{V}_\delta 
\big(M_{>\widehat{r}_\lambda^-}^*
\big)_{\widehat{ r}_\lambda^-}
\bigcup
\widetilde{\Omega}_{>\widehat{r}_\lambda^-}^*.
\]

To establish {\bf (iv)} of Proposition~5.3 at $( \cdot )_{ > \widehat{
r}_\lambda^-}$, it suffices to show {\bf (iii)}, which is checked to
be equivalent. We observe that, for logical reasons only, a given
region $\widetilde{ \sf R}_{ \widehat{ r}_\lambda^+ }^c$ for $c \neq
c_1$ and $c \neq c_2$ can:

\begin{itemize}

\smallskip\item[$\bullet$]
be disjoint from $\widetilde{ \sf R}_{ \widehat{ r}_\lambda^+ }^{
c_1}$ and also disjoint from $\widetilde{ \sf R}_{ \widehat{
r}_\lambda^+ }^{ c_2}$;

\smallskip\item[$\bullet$]
be contained in $\widetilde{ \sf R}_{ \widehat{ r}_\lambda^+ }^{ c_1}$
or (exclusive ``or'') in $\widetilde{ \sf R}_{ \widehat{ r}_\lambda^+
}^{ c_2}$;

\smallskip\item[$\bullet$]
contain $\widetilde{ \sf R}_{ \widehat{ r}_\lambda^+ }^{ c_1}$ or
(inclusive ``or'') $\widetilde{ \sf R}_{ \widehat{ r}_\lambda^+ }^{
c_2}$.

\end{itemize}\smallskip

But we claim that in the latter case, $\widetilde{ \sf R}_{ \widehat{
r}_\lambda^+ }^c$ necessarily contains both regions $\widetilde{ \sf
R}_{ \widehat{ r}_\lambda^+ }^{ c_1}$ and $\widetilde{ \sf R}_{
\widehat{ r}_\lambda^+ }^{ c_2}$. Indeed, otherwise the boundary
${\sf N}_{ \widehat{ r}_\lambda^+ }^c$ of $\widetilde{ \sf R}_{
\widehat{ r}_\lambda^+ }^c$ should separate $\widetilde{ \sf R}_{
\widehat{ r}_\lambda^+ }^{ c_1} \cap {\sf C}_\eta$ from $\widetilde{
\sf R}_{ \widehat{ r}_\lambda^+ }^{ c_2} \cap {\sf C}_\eta$ in the
level set $\big\{ v = \frac{ \eta}{ 2} \big\} \cap {\sf C}_\eta$,
which is impossible since ${\sf N}_{ \widehat{ r}_\lambda^+ }^c \cap
{\sf C}_\eta$ is exactly equal to $\big( { \sf N}_{ \widehat{
r}_\lambda^+ }^{ c_1, 1} \cap {\sf C}_\eta \big) \bigcup \big( { \sf
N}_{ \widehat{ r}_\lambda^+ }^{ c_2, 1} \cap {\sf C}_\eta \big)$, not
more.

It follows in all cases that ${\sf N}_{ \widehat{ r}_\lambda^+ }^c =
\partial \widetilde{ \sf R}_{ \widehat{ r}_\lambda^+ }^c$ is disjoint
from ${\sf C}_\eta$, hence it lies in $\big\{ \widehat{ r}_\lambda^-
\leqslant \vert \! \vert z \vert \! \vert \leqslant \widehat{
r}_\lambda^+ \big\} \big \backslash {\sf C}_\eta$. Consequently, the
regular flow of $\frac{ \nabla \, r_M}{ \vert \! \vert \nabla r_M
\vert \! \vert}$ on
\[
\left[
M\cap
\big\{ 
\widehat{r}_\lambda^-
\leqslant 
\vert\!\vert z 
\vert\!\vert 
\leqslant 
\widehat{r}_\lambda^+ 
\big\} 
\right] 
\big
\backslash {\sf C}_\eta
\]
pushes down regularly ${\sf N}_{ \widehat{ r }_\lambda^+ }^c$, as a
uniquely defined compact 2-codimensional ${\sf N}_{ \widehat{ r
}_\lambda^- }^c \subset {\sf S}_{\widehat{ r }_\lambda^+ }^{ 2n - 1}
$, disjointly from the newly created merged boundary ${\sf N}_{
\widehat{ r}_\lambda^- }^* = \partial \widetilde{ \Omega}_{ >
\widehat{ r}_\lambda^-}^* \subset {\sf S}_{ \widehat{ r}_\lambda^- }^{
2n - 1}$. This information suffices now to check that {\bf (iii)} and
{\bf (iv)} of Proposition~5.3 are transmitted to $( \cdot )_{ >
\widehat{ r}_\lambda^- }$, just for logical reasons.

\smallskip

The proof of {\sf B} in case $k_\lambda = 1$, subcase {\sf (a)} is
complete. Subcase {\sf (a')} involves only minor differences.

\subsection*{6.3.~Subtracting process}
We now summarize the discussion of subcase {\sf (b)}, focusing only on
topological aspects and dropping the formal considerations about
holomorphic extensions. For an adequate three-dimensional
illustration, think of a smoothly cut cylindrical piece of modelling
clay in which a thin finger drills a hole.

\begin{center}
\input involute-saddle.pstex_t
\end{center}

As in \S5.6, in ${\sf C}_\eta$, there enter exactly two domains
$\widetilde{ \Omega}_{> \widehat{ r }_\lambda^+ }^{ c_i}$, $i=1, 2$,
with $\widetilde{ \Omega }_{ > \widehat{ r }_\lambda^+ }^{ c_1} \subset
\widetilde{ \Omega }_{ > \widehat{ r }_\lambda^+ }^{ c_2}$ by the
induction assumption. Also, there enter two connected regions
$\widetilde{ \sf R}_{ \widetilde{ r}_\lambda^+}^{ c_i, 1} \subset {\sf
S}_{ \widehat{ r}_\lambda^+}^{ 2n - 1}$, $i=1, 2$, with $\widetilde{
\sf R}_{ \widetilde{ r}_\lambda^+ }^{ c_1, 1} \subset \widetilde{ \sf
R}_{ \widetilde{ r}_\lambda^+ }^{ c_2, 1}$. Their boundaries contain
two connected hypersurfaces ${\sf N}_{ \widehat{ r}_\lambda^+}^{ c_i,
1}$ of ${\sf S}_{ \widehat{ r}_\lambda^+}^{ 2n - 1}$, $i = 1, 2$,
which enter ${\sf C}_\eta$ as the two caps of the third pic of
Figure~19.

By descending the interval $( \widehat{ r}_\lambda, \widehat{
r}_\lambda^+ )$ up to $( \cdot )_{ > \widehat{ r}_\lambda }$, we get
two regions $\widetilde{ \sf R}_{ \widehat{ r}_\lambda }^{ c_i, 1}$,
$i = 1, 2$, that touch at $\widehat{ p }_\lambda$, namely the left
cone and the exterior of the right cone in the second pic of
Figure~19.

While descending further to $( \cdot)_{ > \widehat{ r}_\lambda -
\varepsilon }$, with $\varepsilon > 0$ very small, the left cone does
merge with the right (white) cone. Observe that the points of this
(white) cone may be joined continuously to points of the (white) right
cap of the first pic, which by hypothesis lies outside $\widetilde{
\Omega}_{ > \widehat{ r}_\lambda^+ }^{ c_2}$, hence in the same
connected component as the points at infinity. Consequently, we cannot
prolong the left domain $\widetilde{ \Omega}_{ > \widehat{ r }_\lambda^+
}^{ c_1}$ so that its prolongation contains the left cone in the slice
$\{ v = 0 \}$ (third pic), because no admissible prolongation would
enjoy the relative compactness {\bf (i)} of Proposition~5.3. Hence we
have no other choice except to suppress $\widetilde{ \Omega}_{ >
\widehat{ r}_\lambda}^{ c_1}$ when attaining $( \cdot)_{ > \widehat{ r
}_\lambda}$. We then get a new domain $\widetilde{ \Omega}_{ >
\widehat{ r}_\lambda }^*$ defined as $\widetilde{ \Omega}_{ >
\widehat{ r}_\lambda}^{ c_2}$ minus the closure of $\widetilde{
\Omega}_{ > \widehat{ r}_\lambda}^{ c_1}$ (subtraction process),
which is checked to be relatively compact in $\C^n$. This domain then
descends as a uniquely defined domain $\widetilde{ \Omega}_{ >
\widehat{ r}_\lambda^- }^*$ at $( \cdot)_{ > \widehat{ r}_\lambda^-
}$. We also get a corresponding connected region $\widetilde{ \sf
R}_{ \widehat{ r}_\lambda^- }^*$ approximately equal to $\widetilde{
\sf R}_{ \widehat{ r}_\lambda}^{ c_2, 1}$ minus the closure of
$\widetilde{ \sf R}_{ \widehat{ r}_\lambda}^{ c_1, 1}$ whose boundary
contains a connectedd ${\sf N}_{ \widehat{ r}_\lambda^- }^*$ (bottom
right Figure~21), obtained by merging ${\sf N}_{ \widehat{
r}_\lambda}^{ c_1, 1}$ with ${\sf N}_{ \widehat{ r}_\lambda}^{ c_2,
1}$.

\smallskip

The last subcase {\sf (e)} above is topologically similar to what
happens in \S5.15, hence the proof of Proposition~5.3 is complete.
\endproof

\vfill\end{document}